\documentclass[12pt]{article}
\usepackage{amsmath,amsfonts,amssymb,amscd}
\usepackage{graphicx}

\newtheorem{theo}{Theorem}
\newtheorem{lem}{Lemma}[section]
\newtheorem{prop}{Proposition}[section]

\newtheorem{cor}{Corollary}[section]
\newtheorem{rem}{Remark}[section]
\newtheorem{dfn}{Definition}[section]

\makeatletter \@addtoreset{equation}{section} \makeatother

\newcommand{\pr}{\mathop{\rm pr}\nolimits}

\newcommand{\mC}{\mathbb{C}}

\newcommand{\mR}{\mathbb{R}}

\newcommand{\mZ}{\mathbb{Z}}
\newcommand{\mN}{\mathbb{N}}

\newcommand{\bF}{{\bf F}}

\newcommand{\bH}{{\bf H}}

\newcommand{\bM}{{\bf M}}

\newcommand{\bV}{{\bf V}}

\newcommand{\be}{{\bf e}}

\newcommand{\bk}{{\bf k}}
\newcommand{\bl}{{\bf l}}
\newcommand{\bm}{{\bf m}}

\newcommand{\bs}{{\bf s}}
\newcommand{\bu}{{\bf u}}
\newcommand{\bv}{{\bf v}}
\newcommand{\bw}{{\bf w}}

\newcommand{\bz}{{\bf z}}

\newcommand{\calA}{{\cal A}}

\newcommand{\calD}{{\cal D}}
\newcommand{\calF}{{\cal F}}

\newcommand{\calH}{{\cal H}}

\newcommand{\calN}{{\cal N}}

\newcommand{\calP}{{\cal P}}

\newcommand{\calW}{{\cal W}}

\newcommand{\eps}{\varepsilon}
\newcommand{\ph}{\varphi}

\newcommand{\lev}{\operatorname{lev}}

\newcommand{\one}{{\bf{1}}}

\newcommand{\re}{\operatorname{Re}}

\newcommand\qed{{\unskip\nobreak\hfil\penalty50
  \hskip2em\hbox{}\nobreak\hfil\mbox{\rule{1ex}{1ex} \qquad}
    \parfillskip=0pt \finalhyphendemerits=0\par\medskip}}

\begin{document}

\title
{Normalization flow}
\author{ Dmitry Treschev \\
Steklov Mathematical Institute of Russian Academy of Sciences
}
\date{}
\maketitle

\begin{abstract}
We propose a new approach to the theory of normal forms for Hamiltonian systems near a non-resonant elliptic singular point. We consider the space of all Hamiltonian functions with such an equilibrium position at the origin and construct a differential equation in this space. Solutions of this equation move Hamiltonian functions towards their normal forms. Shifts along the flow of this equation correspond to canonical coordinate changes. So, we have a continuous normalization procedure. The formal aspect of the theory presents no difficulties. The analytic aspect and the problems of convergence of series are as usual non-trivial.
\end{abstract}

\section{Introduction}

Following Birkhoff, \cite{Birk} consider the problem of reduction of a Hamiltonian system in a neighborhood of an elliptic fixed point to the normal form. More precisely, we deal with the system
\begin{equation}
\label{ham_sys}
   \dot z = i\partial_{\overline z} \widehat H, \quad
   \dot{\overline z} =  - i\partial_z \widehat H, \qquad
   \widehat H = \widehat H(z,\overline z).
\end{equation}
Here $z = (z_1,\ldots,z_n)$ and $\overline z = (\overline z_1,\ldots,\overline z_n)$ are independent coordinates on $\mC^{2n}$, the Hamiltonian function $\widehat H$ has the form\footnote
{Note that the variable $\overline z_j$ is not necessarily complex conjugate of $z_j$. Following Birkhoff, we use $z$ and $\overline z$ as {\it independent} complex variables. Analogously the notation $\overline k$ has nothing to do with complex conjugacy.}
\begin{eqnarray}
\label{ham_ini}
&\!\!\!\!\!\!\!\!
  \widehat H = \sum_{k,\overline k\in\mZ_+^n} \widehat H_{k,\overline k} z^k\overline z^{\overline k}
             = H_2 + \widehat H_\diamond, \qquad
  \widehat H_\diamond = O_3(z,\overline z), \quad
  H_2(z,\overline z) = \sum_{j=1}^n \omega_j z_j\overline z_j, & \\
\nonumber
& z^k = z_1^{k_1}\cdots z_n^{k_n}, \quad
  \overline z^{\overline k} = \overline z_1^{\overline k_1}\cdots \overline z_n^{\overline k_n}, \quad
  \mZ_+^n = \{0,1,2,\ldots \}. &
\end{eqnarray}
Here $O_m(z,\overline z)$ is a short form for $O(|z|^m + |\overline z|^m)$.

Reality condition for $\widehat H$ is as follows:
\begin{equation}
\label{r-a}
  \overline{\widehat H}_{k,\overline k} = \widehat H_{\overline k, k} \quad
  \mbox{for any $k,\overline k\in\mZ_+^n$}.
\end{equation}
In general we do not assume that $\widehat H$ is real in this sense.

We assume that the frequency vector $\omega = (\omega_1,\ldots,\omega_n)$ is nonresonant:
$$
  \langle q,\omega\rangle \ne 0 \quad\mbox{for any $q\in\mZ^n\setminus\{0\}$},
$$
where $\langle\cdot,\cdot\rangle$ is the standard inner product on $\mR^n$.
It is known that in this case there exists a formal canonical near-identity change of the variables
\begin{equation}
\label{change}
  (z,\overline z)\mapsto (Z,\overline Z) = (z,\overline z) + O_2(z,\overline z), \quad
  d\overline z\wedge dz = d\overline Z\wedge dZ
\end{equation}
such that the new Hamiltonian function takes the normal form:
$$
  H(z,\overline z) = N(Z_1\overline Z_1,\ldots,Z_n\overline Z_n).
$$

The change of variables (\ref{change}) is formal because $Z = Z(z,\overline z)$,
$\overline Z = \overline Z(z,\overline z)$, and $N$ are in general formal power series.
The problem of convergence/divergence of the normalizing transformation under the assumption of analyticity of $\widehat H$ is central in the theory. A separate (harder) problem is convergence/divergence of the normal form. If the normalization converges then the system is locally completely integrable. Various versions of the inverse statement are proved in \cite{Vey,Ito,El}.

Another corollary from convergence of the normalization is Lyapunov stability of the equilibrium position. The papers \cite{Fayad,Koz} contain examples of real-analytic Hamiltonians (\ref{ham_ini}) such that the origin is Lyapunov unstable in the system (\ref{ham_sys}).

Convergence of the normal form does not imply convergence of the normalization. But it has interesting dynamical consequences: the measure of the set covered by KAM-tori turns out to be noticeably bigger than in the case when the normal form diverges \cite{Krik}.

Probably one should expect that convergence is an exceptional phenomenon in any reasonable sense. At the moment this exceptionality is known in terms of Baire category \cite{Sie} and $\Gamma$-capacity \cite{PM,Krik}. Explicit examples of systems with divergent normal form can be found in \cite{Gong,Yin,Fayad}.

The case of the ``trivial'' normal form
\begin{equation}
\label{triv}
   N(Z_1\overline Z_1,\ldots,Z_n\overline Z_n) = H_2(Z,\overline Z)
\end{equation}
is special. According to \cite{Bryuno,Russ}, see also \cite{Stol}, (\ref{triv}) combined with Diophantine conditions, imposed on $\omega$, imply convergence of the normalization.

The change of variables (\ref{change}) is traditionally constructed as a composition of an infinite sequence of coordinate changes which normalize the Hamiltonian function up to a remainder of a higher and higher degree \cite{Birk}. In further works (see for example \cite{Si-Mo}) the change of coordinates (\ref{change}) is represented as a formal series in $z$ and $\overline z$. Because of the symplecticity condition the corresponding formulas look cumbersome and implicit, but such an approach gives a possibility to control polynomial structure of the normalization in any finite degree w.r.t. the phase variables.
\smallskip

Let $\calF$ be the space of all power series in the variables $z$ and $\overline z$. Sometimes we refer to elements of $\calF$ as functions although they are only formal power series.
In this paper we study a flow on the space\footnote
{More precisely on the subspace $\calF_\diamond\subset\calF$ of series which start from terms of power 3.}
$\calF$. We call it the normalization flow $\phi^\delta$, $\delta\ge 0$. Any shift
$$
  H_2 + \widehat H_\diamond \mapsto H_2 + \phi^\delta(\widehat H_\diamond), \qquad
  \widehat H_\diamond\in\calF
$$
is a transformation of the Hamiltonian function $H_2 + \widehat H_\diamond$ according to a certain (depending on $\delta$ and on the initial Hamiltonian $H_2 + \widehat H_\diamond$) canonical change of variables.

The flow $\phi^\delta$ is determined by a certain ODE in $\calF$:
\begin{equation}
\label{main}
    \partial_\delta H_\diamond
  = - \{\xi H_\diamond, H_2 + H_\diamond\}, \qquad
    H_\diamond|_{\delta=0} = \widehat H_\diamond.
\end{equation}
Here $\{\,,\}$ is the Poisson bracket and $\xi$ is a linear operator on $\calF$. In fact, (\ref{main}) is an initial value problem (IVP) for a differential equation presented in the form of Lax L-A pair. Usually such systems are considered in the theory of integrable systems.

Let $\calN\subset\calF$ be the subspace which consists of series which depend on $z,\overline z$ only through the combinations $z_j\overline z_j$. The space $\calN$ is an invariant manifold with respect to $\phi^\delta$. Any point of $\calN$ is fixed. Any orbit of the flow indefinitely approaches $\calN$ with respect to a certain topology on $\calF$.

In Section \ref{sec:xi} we represent the system (\ref{main}) in the form of an infinite ODE system for the coefficients $H_{k,\overline k}$. A bit more convenient equivalent form of it is the system (\ref{aver4}). Because of a special ``nilpotent'' structure of the system (\ref{aver4}), the existence and uniqueness of a solution for the corresponding IVP for any initial condition turns out to be a simple fact (Section \ref{sec:formal}).

We are particularly interested in the restriction of $\phi^\delta$ to the subspace $\calA\subset\calF$ of analytic Hamiltonian functions. In Section \ref{sec:analy} we prove that for any
$\widehat H_\diamond\in\calA$ the solution $H_\diamond$ also lies in $\calA$ for any $\delta>0$. However the polydisk of analyticity generically shrinks when $\delta$ grows. According to our rough lower estimate its radius is of order $1/\delta$.

In Section \ref{sec:simple} we discuss some special properties of the system (\ref{main}). For example, the set
$$
  \{H_\diamond\in\calF : H_{k,\overline k}=0 \mbox{ if } \langle\omega,\overline k - k\rangle\le M_1
                                            \mbox{ or } \langle\omega,\overline k - k\rangle\ge M_2 \}
$$
is invariant for the flow $\phi^\delta$.

In Section \ref{sec:advanced} after some technical work we present the system (\ref{aver4}) in a more convenient form. In this form unknown functions in the differential equations are the functions $\calH^q=\sum_{\overline k-k=q}\calH_{k,\overline k} z^k\overline z^{\overline k}$, $q\in\mZ^n$. Any function $\calH^q$ equals $z^{k_q}\overline z^{\overline k_q} N^q$. Here
$z^{k_q}\overline z^{\overline k_q}$ is a certain monomial, determined by $q$, and $N^q\in\calN$. We rewrite the system in terms of the functions $N^q$ only. This version of the system (\ref{main}) is the main tool in Sections \ref{sec:expli} and \ref{sec:exa} to obtain some explicit solutions and in Section \ref{sec:conver} to prove a global in $\delta\ge 0$ existence theorem for a certain initial condition $\widehat H_\diamond$.

For any $N\in\calN$ let $\calW_N$ be the corresponding stable manifold. Then $\calF$ is foliated by such manifolds: $\calF = \cup_{N\in\calN} \calW_N$. Unfortunately the invariant manifold $\calN$ is not normally hyperbolic. Because of small divisors the orbits $\phi^\delta(\widehat H_\diamond)$ do not tend to their normal forms uniformly exponentially. However it seems interesting to consider the restriction of $\phi^\delta$ to the manifolds $\calW_N$. Section \ref{sec:asy} contains an attempt to imagine how such a restriction might look like.

This is a motivation to introduce the so-called asymptotic flow $\Psi^\delta$. It is determined by an equation (asymptotic system) analogous to (\ref{main}), where in comparison with (\ref{main}) some terms are dropped. Due to this the ``normal component'' of the Hamiltonian
$\Psi^\delta(\widehat H_\diamond)$ remains constant. In Section \ref{sec:asy_exi} we prove that the asymptotic flow is conjugated to the normalization flow (presented in slightly another form $\Phi^\delta$) by a smooth map $\calF\to\calF$.

In Section \ref{sec:expli} we obtain an explisit solution of the asymptotic system. At this point one might expect that at least for the asymptotic flow analytic initial conditions generate an orbit on which the analyticity domain does not shrink to zero. Unfortunately, this is not true. In Section \ref{sec:exa} we consider an example when the normalization flow (in the form $\Phi^\delta$) and the asymptotic flow have the same orbit. This happens when all Taylor coefficients of the initial condition $\widehat H_\diamond$ with $\langle\omega,\overline k-k\rangle < 0$ vanish. According to results of Section \ref{sec:simple} this property remains valid on the whole orbit. Note that such an example is physically meaningless because this choice of $\widehat H_\diamond$ is incompatible with the reality condition. However in the context of our approach it is important because we prove (Proposition \ref{prop:asym_diver}) that typically\footnote
{Even in the case when the normal form (which is known in this example in advance) is a quadratic polynomial in $\kappa_j = z_j\overline z_j$.}
the analyticity polydisk for the solution $H_\diamond$ shrinks indefinitely as $\delta\to\infty$.

Section \ref{sec:conver} presents a somewhat more positive result. We present an example of an analytic (if necessary, real-analytic) Hamiltonian
$\widehat H_\diamond$ for which both its normal form $\phi^\infty(\widehat H_\diamond)$ and the transformation to the normal form are analytic. Such a function $\widehat H_\diamond$ satisfies the following condition:
$$
  \widehat H_{k,\overline k} \ne 0\quad\mbox{only for}\quad \overline k - k\in\{0,q,-q\},
$$
where $q\in\mZ^n$ is a nonzero vector. Hence, ``majority'' of coefficients
$\widehat H_{k,\overline k}$ in this case vanish, although expansion of $\widehat H_\diamond$ still contains an infinite number  of nonzero terms.

Analytic convergence to normal form in this situation is not very surprising. Indeed, any quadratic form
$$
  Q = \sum a_j z_j\overline z_j , \qquad
  \sum a_j q_j = 0
$$
commutes with $\widehat H_\diamond$: $\{Q,\widehat H_\diamond\}=0$. Therefore the traditional normalization procedure converges by the Ito theorem \cite{Ito}. However we do not know if normalization by the flow $\phi^\delta$ is equivalent from the viewpoint of analytic convergence to the traditional normalization. On the other hand it was interesting for us  to prove a global in time $\delta$ existence of an analytic solution of IVP (\ref{main}) for some nontrivial initial condition $\widehat H_\diamond$. Technically the problem is reduced to a certain nonlocal in time theorem of Cauchy-Kovalevskaya type for the PDE system (\ref{pde_eq1}).

In this paper dealing with problems of local in ``space'' but global in ``time'' analytic convergence of solutions for various initial value problems we use the so-called Majorant principle presented in Section \ref{sec:MP}. It slightly differs from majorant argument traditionally used in relation with the Cauchy-Kovalevskaya theorem. The main difference is that we refuse to regard time as a complex variable and therefore no Taylor expansions in time appear. This version of the majorant method is better adapted to the problem of the construction of nonlocal in time analytic solutions.

Finally note that we discuss a ``continuous'' (unlike step-by-step) approach to the theory of normal forms under the following restrictions:
\begin{itemize}
\item the systems are Hamiltonian,
\item the fixed point is elliptic,
\item the frequencies are nonresonant.
\end{itemize}
Surely our methods work if any of these conditions (or any combination of these conditions) is dropped.
\smallskip

{\bf Acknowledgements}. The author thanks Sergey Bolotin, Sergey Suetin and Oleg Zubelevich for valuable discussions.

\section{Basic construction}

For any $q = (q_1,\ldots,q_n)\in\mZ^n$ let $|q|=|q_1|+\cdots+|q_n|$ be its $l^1$-norm. For any
$k,\overline k\in\mZ_+^n$ we put
\begin{eqnarray*}
& \bk=(k,\overline k)\in\mZ_+^{2n},\quad
  \bk' = \overline k - k\in\mZ^n , \quad
  \bk^* = (\overline k,k), \quad
  |\bk| = |k| + |\overline k|, & \\
& \bz=(z,\overline z),\quad
  \bz^\bk = z^k\overline z^{\overline k},\quad
  z,\overline z\in\mC^n . &
\end{eqnarray*}
Then $\bk^*-\bk=(\bk',-\bk')$. In particular $\bk'=0$ if and only if $\bk=\bk^*$.

We will use the notation
$$
  \mZ_\diamond^{2n} = \{\bk\in\mZ_+^{2n} : |\bk|\ge 3\}.
$$

\subsection{Spaces}

Let $\calF$ be the vector space of all series
\begin{equation}
\label{H_*incalF}
  H = \sum_{\bk\in\mZ_+^{2n}} H_{\bk} \bz^\bk, \qquad
  H_{\bk}\in\mC.
\end{equation}
The series (\ref{H_*incalF}) are assumed to be formal i.e., there is no restriction on the values of the coefficients $H_{\bk}$. So, $\calF$ coincides with the ring
$\mC[[z_1,\ldots,z_n,\overline z_1,\ldots,\overline z_n]]$.

Below we use on $\calF$ the product topology i.e., a sequence $H^{(1)},H^{(2)},\ldots\in\calF$ is said to be convergent if for any $\bk\in\mZ_+^{2n}$ the sequence of coefficients
$H^{(1)}_{\bk},H^{(2)}_{\bk},\ldots$ converges.

For any $H$ satisfying (\ref{H_*incalF}) we define $p_{\bk}(H) = H_{\bk}$. So $p_{\bk} : \calF\to\mC$ is a canonical projection corresponding to $\bk\in\mZ_+^{2n}$.
Suppose $F\in\calF$ depends on a parameter $\delta\in I$, where $I\subset\mR$ is an interval. In other words, we consider a map
$$
  f : I\to\calF, \quad I\ni\delta\mapsto F(\cdot,\delta).
$$
We say that $F$ is smooth in $\delta$ if all the maps $p_\bk\circ f$ are smooth.


Let $\calF_r\subset\calF$ be the space of ``real'' series:
$$
   \calF_r
 = \{H\in\calF : \overline H_{\bk} = H_{\bk^*} \;
                                    \mbox{ for any $\bk\in\mZ_+^n$}\}.
$$

We define $\calA\subset\calF$ as the space of analytic series:
$$
   \calA
 = \{H\in\calF : \mbox{ there exist $c,a$ such that $|H_{\bk}| \le c e^{a|\bk|}$
                          for any $\bk\in\mZ_+^{2n}$}\}.
$$

The product topology may be restricted from $\calF$ to $\calA$, but, being a scale of Banach spaces, $\calA$ may be endowed with a more natural topology. We have
$\calA = \cup_{0<\rho}\calA^\rho$, where $\calA^\rho$ is a Banach space with the norm $\|\cdot\|_\rho$,
\begin{equation}
\label{Drho}
  \|H\|_\rho = \sup_{\bz\in D_\rho} |H(\bz)|, \quad
  D_\rho = \{\bz\in\mC^{2n} : |z_j|<\rho,\; |\overline z_j|<\rho,\; j=1,\ldots,n\}.
\end{equation}
Then we have:
$\calA^{\rho'}\subset\calA^\rho$, $\|\cdot\|_\rho \le \|\cdot\|_{\rho'}$ for any
$0 < \rho < \rho'$.

\begin{lem}
\label{lem:Hkk}
If $\|H\|_\rho \le c$ then
\begin{equation}
\label{Cauchy}
  |H_{\bk}| \le c\rho^{-|\bk|}.
\end{equation}
\end{lem}

The proof follows from the Cauchy formula: for any positive $\rho_0<\rho$ and any $H\in\calA^\rho$
$$
    H_{\bk}
  = \frac1{(2\pi i)^{2n}} \oint dz_1 \ldots \oint dz_n \oint d\overline z_1 \ldots
                  \oint \frac{H(\bz)\, d\overline z_n}{\bz^{\bk+\one}},
$$
where $\one=(1,\ldots,1)\in\mZ_+^{2n}$ and the integration is performed along the circles
$$
  |z_1|=\rho_0,\ldots,\quad
  |z_n|=\rho_0,\quad
  |\overline z_1|=\rho_0,\ldots,\quad
  |\overline z_n|=\rho_0.
$$
This implies $|H_\bk|\le c\rho_0^{-|\bk|}$. Since $\rho_0<\rho$ is arbitrary, we obtain (\ref{Cauchy}). \qed

\begin{cor}
\label{cor:stronger}
Topology determined on $\calA^\rho$ by the norm $\|\cdot\|_\rho$ is stronger than the product topology induced from $\calF$ i.e., if a sequence $\{H^{(j)}\}$ converges in $(\calA^\rho, \|\cdot\|_\rho)$ then it converges in $\calF$.
\end{cor}

Let $\calN\subset\calF$ be the space of ``normal forms'':
$$
   \calN
 = \{H\in\calF : H_{\bk} \ne 0 \; \mbox{ implies } \bk'=0\}.
$$

We define the subspace (the subring) $\calF_\diamond\subset\calF$ of the series
\begin{equation}
\label{H_diamondincalF}
  H_\diamond = \sum_{\bk\in\mZ_\diamond^{2n}} H_{\bk} \bz^\bk, \qquad
  H_{\bk}\in\mC.
\end{equation}
By definition we also have: $\calN_\diamond = \calN\cap\calF_\diamond$.

\subsection{Continuous averaging}

Now we construct a flow $\phi^\delta$, $\delta\ge 0$ in $\calF_\diamond$ which asymptotically (as $\delta\to+\infty$) reduces any Hamiltonian $H_2 + \widehat H_\diamond$, $\widehat H_\diamond\in\calF_\diamond$ to the normal form $H_2 + N_\diamond$, $N_\diamond=\phi^{+\infty}(\widehat H_\diamond)\in\calN_\diamond$. Construction of the normalization flow is based on the method of continuous averaging \cite{TZ}.

First we explain the general idea of the continuous averaging. Consider the change of variables in the form of a shift
\begin{equation}
\label{shift}
  \bz = (z,\overline z) \mapsto \bz_\delta = (z_\delta,\overline z_\delta), \qquad
  (Z,\overline Z) = (z_\Delta,\overline z_\Delta)
\end{equation}
along solutions of the Hamiltonian system\footnote{
The so-called, Lie method.} with Hamiltonian $F=F(\bz,\delta)=O_3(\bz)$ and independent variable $\delta$:
\begin{equation}
\label{ham=F}
  z' = i\partial_{\overline z} F, \quad
  \overline z' = -i\partial_z F, \qquad
  (\cdot)' = d/d\delta.
\end{equation}
Suppose the function $H_2+\widehat H_\diamond$ expressed in the variables $\bz_\delta$ takes the form $H_2+H_\diamond$:
\begin{equation}
\label{H=H}
    H_2(\bz) + \widehat H_\diamond(\bz)
  = H_2(\bz_\delta) + H_\diamond(\bz_\delta,\delta).
\end{equation}

Let $\{\,,\,\}$ be the Poisson bracket on $\calF$:
$$
     \{F,G\}
  = i \sum_{j=1}^n \big( \partial_{\overline z_j}F \partial_{z_j}G
                        - \partial_{z_j}F \partial_{\overline z_j}G \big).
$$
Differentiating equation (\ref{H=H}) in $\delta$, we obtain:
$$
    \partial_\delta H_\diamond
  = -\{F,H_2+H_\diamond\}, \qquad
    H_\diamond|_{\delta=0} = \widehat H_\diamond.
$$

The main idea of the continuous averaging is to take $F$ in the form $\xi H_\diamond$, where $\xi$ is a certain operator\footnote
{This idea is similar to the Moser's idea of homotopy method \cite{Moser}.}
, depending on the problem we deal with. Then we obtain an initial value problem in $\calF_\diamond$:
\begin{equation}
\label{aver0}
  \partial_\delta H_\diamond = - \{\xi H_\diamond,H_2 + H_\diamond\}, \qquad
  H_\diamond|_{\delta=0} = \widehat H_\diamond.
\end{equation}

\subsection{The operator $\xi$}
\label{sec:xi}

We put
$$
   \sigma_q = \mbox{sign\,}(\langle\omega,q\rangle), \quad
   \omega_q = |\langle\omega,q\rangle|, \qquad
   q\in\mZ^n.
$$
For any $H_\diamond$ satisfying (\ref{H_diamondincalF}) we define
\begin{equation}
\label{xi}
 \xi H_\diamond = -i \sum_{|\bk|\ge 3}
                  \sigma_{\bk'} H_{\bk} \bz^\bk
         = i (H^- - H^+), \qquad
   H^\pm = \sum_{\pm\sigma_{\bk'}>0} H_{\bk} \bz^\bk.
\end{equation}
We also put
$$
  H^0 = \sum_{\bk'=0,\,|\bk|\ge 4} H_{\bk} \bz^\bk \in\calN_\diamond
$$
so that $H_\diamond = H^0 + H^- + H^+$.

We obtain a more detailed form of the differential equation (\ref{aver0}):
\begin{eqnarray}
\label{aver1}
&   \partial_\delta H_\diamond
  = v_0(H_\diamond) + v_1(H_\diamond) + v_2(H_\diamond), & \\
\nonumber
&   v_0
  = -i\{H^- - H^+,H_2\}, \quad
    v_1
  = - i\{H^- - H^+,H^0\}, \quad
    v_2
  = - 2i \{H^-,H^+\}.
\end{eqnarray}
Informal explanation for the choice (\ref{xi}) of the operator $\xi$ is as follows. Removing nonlinear terms $v_1+v_2$ in (\ref{aver1}), we obtain the equation
$\partial_\delta H_\diamond = v_0(H_\diamond)$ or, in more detail
\begin{equation}
\label{aver_trun}
    \partial_\delta H_{\bk}
  = -\omega_{\bk'} H_{\bk}, \qquad
    H_{\bk}|_{\delta=0} = \widehat H_{\bk}
\end{equation}
which can be easily solved:
$$
  H_{\bk} = e^{-\omega_{\bk'}\delta} \widehat H_{\bk} .
$$
Thus when $\delta\to +\infty$, solution $H_\diamond$ of the truncated problem (\ref{aver_trun}) tends to $H^0\in\calN_\diamond$.

Let $e_j=(0,\ldots,1,\ldots,0)$ be the $j$-th unit vector in $\mZ_+^n$ and let $\be_j=(e_j,e_j)\in\mZ_+^{2n}$. Equation (\ref{aver1}) written for each Taylor coefficient $H_{\bk}$ has the form
\begin{eqnarray}
\label{aver3}
     \partial_\delta H_{\bk}
 &=& -\omega_{\bk'} H_{\bk} + v_{1,\bk}(H_\diamond) + v_{2,\bk}(H_\diamond),\qquad
     H_\bk|_{\delta=0} = \widehat H_\bk, \\
\label{v0}
     v_{1,\bk}
 &=& - \sum_{j=1}^n \sum_{|l|\ge 2} \sigma_{\bk'}
                   (\overline k_j - k_j) l_j H_{\bk+\be_j-(l,l)} H_{l,l}, \\
\nonumber
     v_{2,\bk}
 &=& 2 \sum_{j=1}^n \sum_{\sigma_{\bl'} < 0 < \sigma_{\bm'},\,
                                 \bl+\bm-\bk = \be_j}
         (\overline l_j m_j - l_j\overline m_j) H_{\bl} H_{\bm}.
\end{eqnarray}

By using the change of variables
$$
  H_{\bk} = \calH_{\bk} e^{-\omega_{\bk'} \delta},
$$
we reduce the equations (\ref{aver3}) to the form
\begin{eqnarray}
\label{aver4}
\!\!\!\!\!\!\!
     \partial_\delta\calH_{\bk}
 &=&  v_{1,\bk}(\calH_\diamond) + \bv_{2,\bk}(\calH_\diamond), \qquad
     \calH_{\bk}|_{\delta=0} = \widehat H_{\bk}, \\
\label{bv2}\!\!\!\!\!\!\!
     \bv_{2,\bk}
 &=& 2 \sum_{j=1}^n \sum_{\sigma_{\bl'} < 0 < \sigma_{\bm'},\,
                                 \bl+\bm-\bk = \be_j}
     \!  (\overline l_j m_j - l_j\overline m_j)\calH_{\bl}\calH_{\bm}
          e^{- \omega_{\bl',\bm'} \delta}, \\
\nonumber\!\!\!\!\!\!\!
     \omega_{\bl',\bm'}
 &=& \omega_{\bl'}+\omega_{\bm'}-\omega_{\bl'+\bm'}
\end{eqnarray}

\begin{rem}

1. The functions $v_{1,\bk}$ and $v_{2,\bk}$ are quadratic polynomials in the variables $H_{\bk}$. The functions $\bv_{2,\bk}$ are quadratic polynomials in $\calH_{\bk}$ with coefficients, depending on $\delta$.

2. The polynomial $v_{1,\bk}$ vanishes for any $\bk=\bk^*$.

3. The polynomials $v_{1,\bk}, v_{2,\bk}$ and $\bv_{2,\bk}$ do not depend on the variables $H_\bs,\calH_{\bs}$, for any $\bs\in\mZ_\diamond^{2n}$ such that $|\bs| > |\bk|-2$.

4. The polynomials $v_{2,\bk}$ and $\bv_{2,\bk}$ do not depend on the variables $H_\bs,\calH_{\bs}$ for any $\bs\in\mZ_\diamond^{2n}$ such that
$0\le\sigma_{\bk'} \langle\omega,\bs'\rangle\le\omega_{\bk'}$.

5. For any $\sigma_{\bl'} < 0 < \sigma_{\bm'}$
$$
    \omega_{\bl',\bm'}
  = \left\{\begin{array}{ccc}
           2\omega_{\bl'} &\mbox{ if }& \sigma_{\bk'} > 0, \\
           2\omega_{\bm'} &\mbox{ if }& \sigma_{\bk'} < 0.
           \end{array}
    \right.
$$
\end{rem}

\section{Normalization flow exists}

\subsection{Formal aspect}
\label{sec:formal}

Let $I\subset\mR$ be an interval containing the point $0$. We say that the curve $\gamma:I\to\calF_\diamond$ is a solution of the system (\ref{aver0}) on the interval $I$ if $\gamma(0)=\widehat H_\diamond$ and the coefficients
$H_{\bk}(\delta)=\pr_{\bk}\gamma(\delta)$ satisfy (\ref{aver3}).

If the solution $\gamma : I\to\calF_\diamond$ exists and is unique then for any $\delta\in I$ we put $\gamma(\delta)=\phi^\delta(\widehat H_\diamond)$.

\begin{dfn}
\label{dfn:nil}
The ordinary differential equation
$$
  F' = \Phi(F), \qquad
  F = \sum_{|\bk|\ge 3} F_\bk \bz^\bk
$$
on $\calF$ has a nilpotent form if for any $\bk\in\mZ_+^{2n}$ we have: $F'_\bk=\Phi_\bk$, where $\Phi_\bk$ is a function of the variables $F_\bm$, $|\bm|<|\bk|$.
\end{dfn}

The system (\ref{aver4}) has a nilpotent form. In particular,
$\partial_\delta\calH_{\bk} = 0$ for $|\bk|=3$.

\begin{theo}
\label{theo:phi}
For any $\widehat H_\diamond\in\calF_\diamond$ and any $\delta\ge 0$ the element $H_\diamond(\cdot,\delta)=\phi^\delta(\widehat H_\diamond)$ is well-defined. For any $\bk\in\mZ_\diamond^{2n}$ the function $\calH_\bk(\delta) = H_\bk(\delta) e^{\omega_{\bk'}\delta}$ equals
\begin{equation}
\label{Pk}
  \calH_\bk(\delta) = \widehat H_\bk + P_\bk(\widehat H_\diamond,\delta),
\end{equation}
where $P_\bk$ is a polynomial in $\widehat H_\bm$, $|\bm|<|\bk|$ with coefficients in the form of (finite) linear combinations of terms $\delta^s e^{-\nu\delta}$, $s\in\mZ_+$. Here $\nu\ge 0$ and moreover, for $\bk'=0$ the equation $\nu=0$ implies $s=0$.
\end{theo}

{\it Proof of Theorem \ref{theo:phi}}. Suppose equations (\ref{Pk}) hold for all vectors $\bk\in\mZ_\diamond^{2n}$ with $|\bk|<K$. Take any $\bk\in\mZ^{2n}_\diamond$ with $|\bk|=K$. By (\ref{aver4})
$$
  \calH_\bk=\widehat H_\bk + I_1 + I_2, \qquad
  I_1 = \int_0^\delta v_{1,\bk}(\calH_\diamond(\lambda))\, d\lambda, \quad
  I_2 = \int_0^\delta
          \bv_{2,\bk}(\calH_\diamond(\lambda),\lambda)\, d\lambda.
$$
By using (\ref{v0}), (\ref{bv2}), and the induction assumption, we obtain (\ref{Pk}).
\qed

\begin{cor}
The limit $\lim_{\delta\to +\infty} \phi^\delta(\widehat H_\diamond)$ exists and lies in $\calN_\diamond$.
\end{cor}

Indeed, by Theorem \ref{theo:phi} for any $\bk\in\mZ_\diamond^{2n}$ there exists the limit
$$
  \lim_{\delta\to +\infty} H_{\bk}(\delta) = H_{\bk}(+\infty).
$$
The convergence is exponential and $H_{\bk}(+\infty)$ vanishes for any $\bk\ne\bk^*$. This is equivalent to the existence of
$\lim_{\delta\to +\infty}\phi^\delta(\widehat H_\diamond)=\phi^{+\infty}(\widehat H_\diamond)$
in the product topology on $\calF_\diamond$ together with the statement
$\phi^{+\infty}(\widehat H_\diamond)\in\calN_\diamond$.

\begin{theo}
\label{theo:real}
Suppose $\widehat H_\diamond\in\calF_r\cap\calF_\diamond$. Then
$\phi^\delta(\widehat H_\diamond)\in\calF_r\cap\calF_\diamond$ for any $\delta\ge 0$.
\end{theo}

{\it Proof}. The required identities $\overline H_{\bk}(\delta) = H_{\bk^*}(\delta)$, $\delta\ge 0$ can be proved by induction in $|\bk|$ by using the equation (\ref{aver4}). \qed

\subsection{Analytic aspect}
\label{sec:analy}

\begin{theo}
\label{theo:analytic}
Suppose $\widehat H_\diamond\in\calA^\rho\cap\calF_\diamond$. Then
$\phi^\delta(\widehat H_\diamond)\in\calA^{g(\delta)}\cap\calF_\diamond$ for any $\delta\ge 0$, where for some constants $A,B>0$
$$
  g(\delta) \ge \frac A{1 + B\delta}.
$$
\end{theo}

{\it Proof}. We use the majorant method. We remind definitions and basic facts concerning majorants in Section \ref{sec:maj}. For some positive $a,b$
$$
      \widehat H_\diamond
  \ll f(\zeta)
   =  \sum_{\bk\in\mZ_\diamond^{2n}} \widehat\bH_{\bk} \bz^\bk, \qquad
  \zeta = \sum_{j=1}^n (z_j + \overline z_j) ,
$$
where the function $f(\zeta)=O_3(\zeta)$, a majorant for the initial condition, will be chosen a bit later.

Consider together with (\ref{aver4}) the following majorant system
\begin{eqnarray}
\label{aver_maj}
\!\!\!\!\!\!
     \partial_\delta\bH_{\bk}
 &=&  \bV_{\bk}, \qquad
     \bH_{\bk}|_{\delta=0} = \widehat\bH_{\bk}, \\
\label{bV}
\!\!\!\!\!\!
     \bV_{\bk}
 &=& 2 \sum_{j=1}^n \sum_{\bl+\bm-\bk = \be_j}
         (\overline l_j m_j + l_j\overline m_j)\bH_{\bl}\bH_{\bm}.
\end{eqnarray}
To obtain equation (\ref{bV}), we have replaced in $v_{1,\bk}$ and $\bv_{2,\bk}$ minuses by pluses, dropped the exponential multipliers, and added some new positive (for positive $\bH_{\bl}$ and $\bH_{\bm}$) terms. The main property of the system (\ref{aver_maj})--(\ref{bV}) is such that if for any $\bk\in\mZ_\diamond^{2n}$ we have
$|\widehat H_{\bk}|\le\widehat\bH_{\bk}$ then for any $\delta\ge 0$ the inequality $|H_\bk|\le\bH_\bk$ holds i.e., conditions (a) and (b) from Definition \ref{dfn:maj} hold.

The system (\ref{aver_maj})--(\ref{bV}) has a nilpotent form. Therefore by Theorem \ref{theo:maj_nil} we may use Majorant principle: if $\bH$ is a solution of (\ref{aver_maj})--(\ref{bV}) then (\ref{aver4}) has a solution $\calH$ and $\calH\ll\bH$ for any $\delta\ge 0$.

The system (\ref{aver_maj})--(\ref{bV}) may be written in a shorter form:
\begin{equation}
\label{maj_shorter}
  \partial_\delta\bH = 4\sum_{j=1}^n \partial_{z_j}\bH\, \partial_{\overline z_j}\bH, \qquad
  \bH|_{\delta=0} = f(\zeta).
\end{equation}
Since the initial condition $\bH|_{\delta=0}$ depends on $z$ and $\overline z$ only through the variable $\zeta$, we may look for a solution of (\ref{maj_shorter}) in the form
$\bH(z,\overline z,\delta)=F(\zeta,\delta)$. The function $F$ satisfies the equation
$$
  \partial_\delta F = 4n (\partial_\zeta F)^2, \qquad
  F|_{\delta=0} = f(\zeta).
$$

The function $G=\partial_\zeta F$ satisfies the inviscid Burgers' equation
\begin{equation}
\label{Burgers}
  \partial_\delta G = 8n G \partial_\zeta G, \qquad
  G|_{\delta=0} = \partial_\zeta f(\zeta) = O_2(\zeta).
\end{equation}

By using the method of characteristics we obtain that the function $G=G(x,t)$ which solves (\ref{Burgers}), satisfies the equation
\begin{equation}
\label{Bur_sol}
  G = f'(\zeta + 8n\delta G).
\end{equation}

\begin{lem}
\label{lem:Burgers}
Suppose the function $f'$ is analytic at zero. Then for any $t\ge 0$ the function $G$, solving (\ref{Bur_sol}), is also analytic at zero. The corresponding radius of analyticity is of order $1/\delta$.
\end{lem}

{\it Proof}. By Lemma \ref{lem:maj1} we can take $f'(\zeta) = a\zeta^2 / (b-\zeta)$. Then putting $\tau=8n\delta$, we obtain from (\ref{Bur_sol})
$$
  G = \frac{a(\zeta+\tau G)^2}{b - \zeta - \tau G}.
$$
This is a quadratic equation w.r.t. $G$. The solution is
$$
  G = \frac{2a\zeta^2}{b-\zeta-2a\tau\zeta + \sqrt{(b-\zeta-2a\tau\zeta)^2 - 4a\tau\zeta^2(1+a\tau)}}.
$$
It is analytic for
$$
  |\zeta| < \frac{b}{1 + 2a\tau + 2\sqrt{a\tau(1+a\tau)}} .
$$
\qed

\begin{cor}
\label{cor:ra}
Suppose $\widehat H_\diamond\in\calA\cap\calF_r$. Then
$\phi^\delta(\widehat H_\diamond)\in\calA\cap\calF_r$ for any $\delta\ge 0$.
\end{cor}

It is important to note that for a typical $\widehat H_\diamond\in\calA$ one should expect that
$\phi^{+\infty}(\widehat H_\diamond)\in\calN_\diamond$ does not belong to $\calA$, \cite{Krik}.

\section{Simple properties of $\phi^\delta$}
\label{sec:simple}

\subsection{Strips and balls}

We define the strips $S_{M_1,M_2}$ and balls $B_M$:
\begin{eqnarray}
\label{strip}
     S_{M_1,M_2}
 &=& \{\bk\in\mZ_\diamond^{2n} : M_1\le \langle\omega,\bk'\rangle \le M_2\}, \\
\nonumber
     B_M
 &=& \{\bk\in\mZ_\diamond^{2n} : |\bk|\le M \}, \qquad   M\ge 0.
\end{eqnarray}
Sometimes in this notation we use $M_2=+\infty$. Then we take in (\ref{strip})
$\langle\omega,\bk'\rangle < +\infty$. Analogously for $M_1=-\infty$.

We have the following subspaces in $\calF$:
\begin{eqnarray*}
     \calF_{S_{M_1,M_2}}
 &=& \{H_\diamond\in\calF_\diamond : H_\bk=0
                 \mbox{ for any } \bk\in\mZ_\diamond^{2n}\setminus S_{M_1,M_2}\}, \\
     \calF_{B_M}
 &=& \{H_\diamond\in\calF_\diamond : H_\bk=0
                 \mbox{ for any } \bk\in\mZ_\diamond^{2n}\setminus B_M \}.
\end{eqnarray*}

\begin{theo}
\label{theo:strip}
For any $M_1\le M_2$ and $M\ge 0$  the sets $\calF_{S_{M_1,M_2}}$ and $\calF_{B_M}$ are $\phi^\delta$-invariant.
\end{theo}

{\it Proof}. Take any $\bk\in\mZ_\diamond^{2n}\setminus S_{N_1,N_2}$. For definiteness we assume that $\langle\omega,\bk'\rangle > N_2$. Consider the corresponding equation (\ref{aver4}).
In the term $v_{1,\bk}$ (see (\ref{v0})) we have:
$$
  \bk+\be_j-(l,l)\in\mZ_\diamond^{2n}\setminus S_{N_1,N_2}.
$$
The conditions of summation in $\bv_{2,\bk}$
$$
  \sigma_{\bl'} < 0 < \sigma_{\bm'}, \quad
  \bl+\bm-\bk = \be_j
$$
imply $\bk' = \bm' + \bl'$. Therefore
$$
    \langle\omega,\bm'\rangle
  = \langle\omega,\bk'\rangle - \langle\omega,\bl'\rangle
  > \langle\omega,\bk'\rangle
  > N_2.
$$
This implies that $\bm\in\mZ_\diamond^{2n}\setminus S_{N_1,N_2}$.

Hence (by using induction in $|\bk|$), we obtain
$v_{1,\bk}=\bv_{2,\bk}=0$.

In the case of the set $\calF_{B_M}$ the argument is analogous: for any $\bk\in B_M$ we see that any term in $v_{1,\bk}$ and any term in $\bv_{2,\bk}$ contains a multiplier $H_\bm$, $\bm\in B_M$.    \qed

\subsection{Discrete symmetries}

Consider the following two antisymplectic involutions:
$$
  I^\pm(\mC^{2n,0})\hookleftarrow, \qquad
  \bz = (z,\overline z) \mapsto I^\pm(\bz) = \pm(\overline z,z).
$$
For any involution $I^\sigma$, $\sigma\in\{+,-\}$ and any $H\in\calF$ we have:
$$
  H^0\circ I^\sigma = (H\circ I^\sigma)^0,\quad
  H^+\circ I^\sigma=(H\circ I^\sigma)^-,\quad \mbox{and}\quad
  H^-\circ I^\sigma=(H\circ I^\sigma)^+.
$$
Equivalently, $(\xi H)\circ I^\sigma = - \xi (H\circ I^\sigma)$.
\smallskip

Note that in the ``real'' coordinates $(x,y)$, $\sqrt2\,z=y+ix$, $\sqrt2\,\overline z=y-ix$ the map $I^+$ changes sign at $x$ while $I^-$ changes sign at $y$.
\smallskip

The function $H\in\calF$ is said to be $I^\sigma$-invariant if $H\circ I^\sigma = H$.
\smallskip

In particular, $H$ is $I^+$-invariant if it is even w.r.t. $x$ and $I^-$-invariant if it is even w.r.t. $y$. Any function $F\in\calN$ is both $I^+$-invariant and $I^-$-invariant. For example, this holds for the function $H_2$.

Suppose the Hamiltonian function $\widehat H\in\calF$ is $I^\sigma$-invariant and defined in a neighborhood of the origin. Then the corresponding system (\ref{ham_sys}) is said to be $I^\sigma$-reversible. This means that if $\bz(t)$ is a solution of (\ref{ham_sys}) then $I^\sigma(\bz(-t))$ is also a solution of (\ref{ham_sys}).

\begin{theo}
\label{theo:reversible}
Let $I=I^\sigma$, $\sigma\in\{+,-\}$. Suppose $\widehat H$ is $I$-invariant. Then $H$, the solution of the averaging system (\ref{aver0}), is $I$-invariant for any $\delta>0$.
\end{theo}

{\it Proof}. Take composition of (\ref{aver0}) with $I$. We obtain in the left-hand side
$$
  (\partial_\delta H_\diamond)\circ I = \partial_\delta (H_\diamond\circ I),
$$
and in the right-hand side
$$
    -\{\xi H_\diamond, H_2+H_\diamond\}\circ I
  = \{(\xi H_\diamond)\circ I, H_2\circ I + H_\diamond\circ I\}
  = \{-\xi (H_\diamond\circ I), H_2 + H_\diamond\circ I\}.
$$
This means that if $H_\diamond$ is a solution of (\ref{aver0}) then $H_\diamond\circ I$ is also a solution of (\ref{aver0}). Since the initial condition $\widehat H_\diamond$ is assumed to be $I$-invariant and solution of (\ref{aver0}) is unique, Theorem \ref{theo:reversible} follows. \qed

\begin{cor}
The normalization flow $\phi^\delta$ preserves $I^\sigma$-reversibility.
\end{cor}

\section{Advanced form for $v_1$ and $\bv_2$}
\label{sec:advanced}

\subsection{More notation}

We put
$$
  Q_q = \{\bk\in\mZ_\diamond^{2n} : \bk' = q\}.
$$
In any set $Q_q$ there is a unique ``minimal'' element $\bk_q$. It can be defined equivalently by any of two conditions
\smallskip

{\bf (1)} $\bk_q\in Q_q$, $|\bk_q| = \min_{\bk\in Q_q} |\bk|$,

{\bf (2)} $\bk_q=(k_q,\overline k_q)\in Q_q$, $\min\{k_{q_j},\overline k_{q_j}\} = 0$ for any $j=1,\ldots,n$.
\smallskip

For any $q\in\mZ^n$ we have $\bk_q = \bk_{-q}^*$ and
$$
  Q_q = \{\bk\in\mZ_\diamond^{2n} : \bk = \bk_q+(l,l),\; l\in\mZ_+^n\}.
$$

We define the subspaces
$$
  \calF^q = \{F\in\calF_\diamond : F = \sum_{\bk\in Q_q} F_\bk \bz^\bk, \; F_\bk\in\mC\}.
$$
Then $\calF^0=\calN_\diamond$ and $\calF_\diamond = \oplus_{q\in\mZ^n} \calF^q$.

For any $F^q\in\calF^q$ and $F^p\in\calF^p$ we have:
$F^q F^p, \{F^q,F^p\}\in\calF^{q+p}$ and
\begin{equation}
\label{FNbullet}
  F^q = \bz^{\bk_q} N, \qquad
  N\in\calN.
\end{equation}

For any $q,p\in\mZ^n$ we define $q\vartriangleleft p\in\mZ^n$ by
$$
    (q\vartriangleleft p)_j
  = \left\{\begin{array}{cl}
           0 & \mbox{if}\quad q_j p_j \ge 0\quad\mbox{ or }\quad |p_j|<|q_j|, \\
         q_j & \mbox{if}\quad q_j p_j < 0\quad \mbox{and}\quad |q_j|<|p_j|, \\
       q_j/2 & \mbox{if}\quad q_j = - p_j.
           \end{array}
    \right.
$$
For any $s\in\mZ^n/2$ we put
$$
  [s] = (|s_1|,\ldots,|s_n|)\in\mZ^n_+/2.
$$

\begin{lem}
\label{lem:kqkp}
For any $q,p\in\mZ^n$
$$
  \bk_q + \bk_p = \bk_{q+p} + (l,l),\qquad
  l = [q\vartriangleleft p] + [p\vartriangleleft q].
$$
\end{lem}

{\it Proof}. It is sufficient to note that
$$
    ([q\vartriangleleft p] + [p\vartriangleleft q])_j
  = \left\{\begin{array}{cl}
                    0 & \mbox{if}\quad q_j p_j \ge 0, \\
  \min\{|q_j|,|p_j|\} & \mbox{if}\quad q_j p_j < 0.
           \end{array}
    \right.
$$
\qed

\subsection{Poisson brackets}

We put
$$
  \partial = (\partial_1,\ldots,\partial_n), \quad
  \partial_j = \partial_{\kappa_j}, \quad
  \kappa_j = z_j\overline z_j.
$$
The operators $\partial_j$ act on functions (formal series), lying in $\calN$.

\begin{lem}
\label{lem:zz}
For any $q,p\in\mZ^n$ and any $N\in\calN$
\begin{eqnarray}
\label{z=zkap}
      \bz^{\bk_q+\bk_p}
  &=& \bz^{\bk_{q+p}} \kappa^{[q\vartriangleleft p] + [p\vartriangleleft q]}, \\
\label{zN}
      \{\bz^{\bk_q},N\}
  &=& i\bz^{\bk_q} \langle q,\partial\rangle N, \\
\label{zz}
      \{\bz^{\bk_q},\bz^{\bk_p}\}
  &=& i \bz^{\bk_{q+p}} \langle q - q\vartriangleleft p - p + p\vartriangleleft q,\partial\rangle
                       \kappa^{[q\vartriangleleft p] + [p\vartriangleleft q]}.
\end{eqnarray}
\end{lem}

We prove Lemma \ref{lem:zz} in Section \ref{sec:zz}

\begin{lem}
\label{lem:bracket}
For any $p,q\in\mZ^n$ and $N^p,N^q\in\calN$
$$
     \{\bz^{\bk_q} N^q,\bz^{\bk_p} N^p\}
  =  i\bz^{\bk_{q+p}}
      \Big( N^q\kappa^{[q\vartriangleleft p]}
               \langle q,\partial\rangle N^p\kappa^{[p\vartriangleleft q]}
          - N^p\kappa^{[p\vartriangleleft q]}
               \langle p,\partial\rangle N^q\kappa^{[q\vartriangleleft p]}
      \Big).
$$
\end{lem}

Proof follows from equations (\ref{z=zkap})--(\ref{zz}). Indeed,
$\{\bz^{\bk_q} N^q,\bz^{\bk_p} N^p\}=B_1+B_2+B_3$,
$$
      B_1
   =  \bz^{\bk_p} N^q \{\bz^{\bk_q},N^p\}, \quad
      B_2
   =  - \bz^{\bk_q} N^p \{\bz^{\bk_p},N^q\}, \quad
      B_3
   =  N^q N^p \{\bz^{\bk_q},\bz^{\bk_p}\}.
$$

By (\ref{z=zkap}) and (\ref{zN}) we have:
$$
     B_1
  =  i\bz^{\bk_q+\bk_p} N^q \langle q,\partial\rangle N^p
  =  i\bz^{\bk_{q+p}} \kappa^{[q\vartriangleleft p] + [p\vartriangleleft q]} N^q
      \langle q,\partial\rangle N^p.
$$
Analogously
$    B_2
  =  - i\bz^{\bk_{q+p}} \kappa^{[q\vartriangleleft p] + [p\vartriangleleft q]} N^p
      \langle p,\partial\rangle N^q$.
By (\ref{z=zkap}) and (\ref{zz})
$$
     B_3
  =  i\bz^{\bk_{q+p}} N^q N^p
        \kappa^{[q\vartriangleleft p]}\langle q,\partial\rangle \kappa^{[p\vartriangleleft q]}
       - \kappa^{[p\vartriangleleft q]}\langle p,\partial\rangle \kappa^{[q\vartriangleleft p]}.
$$
Lemma \ref{lem:bracket} follows from these computations. \qed

\subsection{Normalizing system revisited}

The equation
$$
     \bv_2(\calH_\diamond,\delta)
  =  -2i \sum_{\sigma_r < 0 < \sigma_s}
           \{\calH^r e^{-\omega_r\delta}, \calH^s e^{-\omega_s\delta}\} e^{\omega_{r+s}\delta}
$$
implies that system (\ref{aver4}) can be presented in the form
\begin{equation}
\label{aver5}
    \partial_\delta\calH^q
  = i\sigma_q\{\calH^q,\calH^0\}
   -2i  \sum_{\sigma_r < 0 < \sigma_s,\, r+s=q}
             \{\calH^r,\calH^s\} e^{-\omega_{r,s}\delta}, \qquad
    \calH^q|_{\delta=0} = \widehat\calH^q .
\end{equation}
Note that $\widehat\calH^q=O_{|q|}(\bz)$ if $|q|\ge 3$, $\widehat\calH^q=O_3(\bz)$ if $0<|q|<3$, and
$\widehat\calH^0=O_4(\bz)$.
\medskip

For any $q\in\mZ^n$ we put
$$
  \calH^q = \bz^{\bk_q} N^q, \quad
  \widehat\calH^q = \bz^{\bk_q}\widehat N^q, \qquad
  N^q,\widehat N^q \in\calN.
$$

By Lemma \ref{lem:bracket}
$$
      \{\calH^r,\calH^s\}
  =  i\bz^{\bk_{r+s}}
      \Big( N^r\kappa^{[r\vartriangleleft s]}
                 \langle r,\partial\rangle N^s\kappa^{[s\vartriangleleft r]}
          - N^s\kappa^{[s\vartriangleleft r]}
                 \langle s,\partial\rangle N^r\kappa^{[r\vartriangleleft s]}
      \Big).
$$
Therefore system (\ref{aver5}) takes the form
\begin{eqnarray}
\label{aver6}
      \partial_\delta N^q
  &=& w_1^q + \bw_2^q, \\
\nonumber
      w_1^q
  &=& - \sigma_q N^q \langle q,\partial\rangle N^0, \\
\nonumber
      \bw_2^q
  &=& 2 \sum_{\sigma_r < 0 < \sigma_s,\, r+s=q}
      \Big( N^r\kappa^{r\vartriangleleft s}\langle r,\partial\rangle N^s\kappa^{s\vartriangleleft r}
          - N^s\kappa^{s\vartriangleleft r}\langle s,\partial\rangle N^r\kappa^{r\vartriangleleft s}
      \Big) e^{-\omega_{r,s}\delta}.
\end{eqnarray}

\section{Asymptotic flow}
\label{sec:asy}

We take $G = \sum_{\bk\in\mZ_\diamond^{2n}} G_\bk(\delta) \bz^\bk$ and
$\widehat G = \sum_{\bk\in\mZ_\diamond^{2n}}\widehat G_\bk \bz^\bk$. The system
\begin{equation}
\label{asym}
   \partial_\delta G_\bk = v_{1,\bk}(G), \qquad
   G_\bk|_{\delta=0} = \widehat G_\bk
\end{equation}
is said to be the asymptotic system. It is obtained from (\ref{aver4}) if we drop the term $\bv_{2,\bk}$. It is motivated by our desire to study restriction of the flow $\ph^\delta$ (or the flow $\Phi^\delta$ of the system (\ref{aver4})) to the stable manifolds $\calW_N$ of points $N\in\calN_\diamond$.

\subsection{Existence theorem }
\label{sec:asy_exi}

Let $G=G(\delta)$ be a solution of the asymptotic system (\ref{asym}). For any $\bk=\bk^*$ we have $v_{1,\bk}=0$. Therefore the coefficients $G_\bk$ with $\bk=\bk^*$ remain constant. Hence we may regard the system (\ref{asym}) as linear. In Section \ref{sec:expli} we obtain an explicit solution for it.

\begin{prop}
Any solution of the asymptotic system has the form
\begin{equation}
\label{pk}
  G_\bk = \widehat G_\bk + p_\bk(\widehat G,\delta),
\end{equation}
where $p_\bk$ are polynomials in $\widehat G_\bm$, $|\bm|<|\bk|$. The dependence of $p_\bk$ on $\delta$ is also polynomial and $\deg_\delta p_\bk \le \lev(\bk)$.
\end{prop}

{\it Proof}.  Induction in $|\bk|$.  \qed

The flow $\Psi^\delta$ of the system (\ref{asym}) will be called the asymptotic flow.

\begin{theo}
\label{theo:conj}
There exists a smooth map $\Lambda:\calF_\diamond\to\calF_\diamond$ which conjugates the flow $\Phi^\delta$ of system (\ref{aver4}) with $\Psi^\delta$:
\begin{equation}
\label{conj}
   \Lambda\circ\Phi^\delta = \Psi^\delta\circ\Lambda.
\end{equation}
\end{theo}

{\it Proof}. Let $\calH_\bk(\delta)$ and $G_\bk(\delta)$ be solutions of (\ref{aver4}) and (\ref{asym}) with the initial conditions
$$
  \calH_\bk(0) = \widehat\calH_\bk \quad \mbox{and}\quad
  G_\bk(0) = \widehat G_\bk.
$$
We determine the map
\begin{equation}
\label{Lambda}
  \Lambda:\calF_\diamond\to\calF_\diamond, \qquad
          \widehat\calH_\diamond = \sum\calH_\bk \bz^\bk
  \mapsto \Lambda(\widehat\calH) = \widehat G_\diamond = \sum G_\bk \bz^\bk.
\end{equation}
by the equation
\begin{equation}
\label{H=G}
  \calH_\bk(\delta) = G_\bk(\delta) + O(e^{-\mu_\bk \delta}), \qquad
  \bk\in\mZ_\diamond^{2n}
\end{equation}
as $\delta\to +\infty$, where $\mu_\bk$ are positive.

\begin{lem}
\label{lem:Lambda}
There exists a unique $\Lambda$ satisfying (\ref{Lambda})--(\ref{H=G}). This map has the form
\begin{equation}
\label{calP}
    \widehat G_\bk
  = \widehat H_\bk + \calP_\bk(\widehat H_\diamond,\delta),
\end{equation}
where $\calP_\bk$ are polynomials in $\widehat H_\bm$, $|\bm|<|\bk|$ with polynomial in $\delta$ coefficients.
\end{lem}

Theorem \ref{theo:conj} obviously follows from Lemma \ref{lem:Lambda}.    \qed
\medskip

{\it Proof of Lemma \ref{lem:Lambda}}. We use induction in $|\bk|$. If $|\bk|=3$, we have:
$$
  \calH_\bk(\delta) = \widehat\calH_\bk \quad \mbox{and}\quad
  G_\bk(\delta) = \widehat G_\bk.
$$
In this case we take $\widehat G_\bk = \widehat\calH_\bk$.

Suppose we have detrmined $\Lambda$ for all vectors from $\mZ_\diamond^{2n}$ with norm less than $|\bk|$. By (\ref{Pk}) and (\ref{pk})
$$
  \calH_\bk(\delta) = \widehat H_\bk + P_\bk(\widehat H_\diamond,\delta), \quad
  G_\bk(\delta) = \widehat G_\bk + p_\bk(\widehat G_\diamond,\delta),
$$
where $P_\bk$ and $p_\bk$ are polynomials in $\widehat H_\bm$ and $\widehat G_\bm$, $|\bm|<|\bk|$. Coefficients of the polynomials $P_\bk$ are linear combinations of expressions
$\delta^s e^{-\nu\delta}$ while coefficients of $p_\bk$ are polynomials in $\delta$.

Let $\widetilde P_\bk$ be the polynomial obtained from $P_\bk$ if all terms
$\delta^s e^{-\nu\delta}$ with positive $\nu$ are replaced by zeros. Then equation (\ref{H=G}) takes the form
$$
    \widehat G_\bk
  = \widehat H_\bk + \widetilde P_\bk(\widehat H_\diamond,\delta)
                   - p_\bk(\widehat G_\diamond,\delta).
$$
By the induction assumption $\widehat G_\bm$, $|\bm|<|\bk|$ have been already expressed as polynomials in $\widehat H_\bl$, $|\bl|<|\bm|$. Hence we obtain (\ref{calP}). \qed

\subsection{Explicit solution}
\label{sec:expli}

To obtain explicit solution for the asymptotic system (\ref{asym}) we present the system in another form. We put
$$
  G = \sum_{q\in\mZ^n} G^q, \qquad
  G^q = \bz^{\bk_q} N^q, \quad
  N^q\in\calN.
$$
Then (\ref{asym}) is equivalent to the system
\begin{equation}
\label{asym2}
    \partial_\delta N^q
  = -\sigma_q N^q \langle q,\partial\rangle N^0, \qquad
    N^q|_{\delta=0} = \widehat N^q.
\end{equation}

Equations (\ref{asym2}) may be solved separately. First we note that $N^0=\widehat N^0$ is independent of $\delta$. Hence
$$
  N^q = e^{-\sigma_q\delta\langle q,\partial\rangle N^0} \widehat N^q, \qquad
  q\in\mZ^n.
$$
Equivalently,
\begin{equation}
\label{sol2}
  G^q = e^{-\sigma_q\delta\langle q,\partial\rangle N^0} \widehat G^q, \qquad
  q\in\mZ^n.
\end{equation}

\subsection{Example}
\label{sec:exa}

Consider the system (\ref{aver3}) with initial condition $\widehat H$ such that for any
$\widehat H_\bk\ne 0$ we have $\bk\in S_{0,+\infty}$. This condition is incompatible with reality of the Hamiltonian $\widehat H$. However we believe that this example is instructive. At least, it rids us from some naive expectations.

By Theorem \ref{theo:strip} $\calH_{\bm}(\delta)=0$ for any
$\bm\in\mZ_\diamond^{2n}\setminus S_{0,+\infty}$ and any $\delta\ge 0$. Therefore $\bv_{2,\bk}=0$ and system (\ref{aver4}) coincides with system (\ref{asym}).

We have: $H^0\equiv\widehat H^0$ and by (\ref{sol2}) for any $q\in\mZ^n$ solution of (\ref{aver3}) reads
\begin{equation}
\label{hq}
  H^q = e^{-(\omega_q + h_q)\delta} \widehat H^q , \qquad
  h_q = \sigma_q \langle q,\partial\rangle\widehat H^0.
\end{equation}
Let $S$ be the set of all $q\in\mZ^n$ such that $\widehat H^q\ne 0$. We put
$$
  c_S := \inf_{q\in S} \omega_q / |q|.
$$

\begin{prop}
Suppose $S\subset \{q\in\mZ^n : \langle\omega,q\rangle \ge 0\}$,
$\widehat H_\diamond\in\calA\cap\calF_\diamond$, and $c_S > 0$.

Then $\widehat H_\diamond$ is transformed to the normal form $\widehat H^0$ by an analytic canonical transformation.
\end{prop}

{\it Proof}. For some $\rho_0>0$ we have: $\widehat H_\diamond\in\calA^{\rho_0}$. Since
$\widehat H^0=O_2(\kappa)$, $\kappa=(\kappa_1,\ldots,\kappa_n)$, $\kappa_j=z_j\overline z_j$, we have:
$$
      \sup_{\bz\in D_\rho} \; \sup_{e\in\mR^n,\,\| e\|=1}
           \big| \langle e,\partial\rangle \widehat H^0(\bz) \big|
  \le c\rho^2,  \qquad
      0 < \rho\le\rho_0,
$$
where $D_\rho$ is the polydisk (\ref{Drho}) and $c$ is a positive constant. Then
$$
  \sup_{\bz\in D_\rho} |h_q(\bz)| \le c |q| \rho^2.
$$
If $\rho^2 < c_S / c$ then
\begin{equation}
\label{infRe}
  \inf_{\bz\in D_\rho} \re (\omega_q + h_q(\bz)) \ge (c_S - c\rho^2) |q| > 0.
\end{equation}
Hence by (\ref{hq}) the functions $H^q$ are defined in $D_\rho$ for any $\delta>0$ and tend to zero as $\delta\to +\infty$.

Inequality (\ref{infRe}) means that there are no small divisors. Hence, change of variables $\bz\mapsto\bz_{+\infty}$ (see (\ref{shift})), generated by the system (\ref{ham=F}) with
$F=\xi H_\diamond$, is well-defined in $D_\rho$ for some positive $\rho$. \qed

The condition $c_S>0$ (no small divisors) is very restrictive.
If $c_S=0$, then typically the change $\bz\mapsto\bz_{+\infty}$ is not analytic and exists only in $\calF$.

\begin{prop}
\label{prop:asym_diver}
Suppose
$$
  S\subset \{q\in\mZ^n : \langle\omega,q\rangle \ge 0\}, \quad
  \widehat H_\diamond\in\calA\cap\calF_\diamond, \quad
  \widehat H^0 = \frac12 \langle A\kappa,\kappa \rangle + O_3(\kappa),
$$
where $A$ is a self-adjoint operator. Let $c$ be a positive constant. We assume that for any $\eps>0$ there exists $q\in S$ such that $\omega_q\le \eps |q|$ and $|Aq| > c|q|$.

Then for any  $\rho>0$ the norm $\|\widehat H_\diamond\|_\rho$ is unbounded as a function of $\delta\in[0,+\infty]$.
\end{prop}

{\it Proof}. Take any $\rho>0$. Let $D_\rho$ be the domain (\ref{Drho}) and
$$
  B_\rho = \{\kappa\in\mC^n : |\kappa_1|<\rho^2, \ldots, |\kappa_n|<\rho^2\}
$$
its image with respect to the map $\bz\mapsto\kappa=(z_1\overline z_1,\ldots,z_n\overline z_n)$.

We have: $h_q = \sigma_q \langle Aq,\kappa\rangle + q O_2(\kappa)$. If $\rho$ is sufficiently small then for some $q\in S$ and an open set $U_\rho\in B_\rho$
$$
  - h_q(\kappa) > \frac12 c |q| \rho^2 > 2\omega_q, \qquad
  \kappa\in U_\rho.
$$
The function $\widehat H^q$ does not vanish identically. Hence
$\sup_{\kappa\in U_\rho} |\widehat H^q| > 0$ and
$$
    H^q(\cdot,\delta)
  = e^{-(\omega_q + h_q)\delta} \widehat H^q \to \infty
  \quad\mbox{for any $\kappa\in U_\rho$}.
$$
This means that for some $\mu>0$ we have: $\|\widehat H^q(\cdot,\delta)\|_{\rho-\mu}\to\infty$ as $\delta\to\infty$. By Lemma \ref{lem:sue} $\|H_\diamond(\cdot,\delta)\|_\rho\to\infty$ as $\delta\to+\infty$. \qed

\section{Converging normalization}
\label{sec:conver}

Take any $q\in\mZ^n\setminus\{0\}$, $\sigma_q>0$. Suppose
$$
  \widehat H_\diamond = \widehat H^q + \widehat H^0 + \widehat H^{-q}, \qquad
  \widehat H^{\pm q} = \bz^{\bk_{\pm q}}\widehat N^{\pm q}\in\calF^{\pm q}, \quad
  \widehat H^0 = \widehat N^0\in\calF^0,
$$
where $\widehat N^{\pm q},\widehat N^0\in\calF^0$.
Then system (\ref{aver6}) contains only three nontrivial equations

\begin{equation}
\label{three}
    \partial_\delta N^{\pm q}
  = - N^{\pm q} \langle q,\partial\rangle N^0, \quad
    \partial_\delta N^0
  = -2 \langle q,\partial\rangle (\kappa^{[q]} N^{-q} N^q) e^{-2\omega_q\delta}.
\end{equation}

Recall that $\widehat H_\diamond = O_3(\bz)$. Therefore $\widehat H^{\pm q} = O_3(\bz)$ and
$\widehat H^0 = O_2(\kappa)$.

\begin{theo}
\label{theo:three_eq}
Suppose $\widehat N^{\pm q},\widehat N^0\in\calN\cap\calA$,
$\bz^{\bk_{\pm q}}\widehat N^{\pm q} = O_3(\bz)$, and $\widehat N^0 = O_2(\kappa)$.
Then there exists $\widehat\rho>0$ such that system (\ref{three}) has a solution
$N^{\pm q},N^0\in\calA^{\widehat\rho}$ for all $\delta>0$. Moreover,
$$
  \lim_{\delta\to+\infty} e^{-\omega_q\delta} N^{\pm q} = 0, \quad
  \lim_{\delta\to+\infty} N^0 \in \calN\cap\calA^{\widehat\rho},
$$
where the limits are taken w.r.t. the norm $\|\cdot\|_{\widehat\rho}$, see (\ref{Drho}).
\end{theo}

{\it Proof}. We put $M = \kappa^{[q]} N^{-q} N^q$. Then
\begin{eqnarray}
\label{M'H0'}
&   \partial_\delta M
  = - 2M \langle q,\partial\rangle N^0, \quad
    \partial_\delta N^0
  = -2 \langle q,\partial\rangle M e^{-2\omega_q \delta}, & \\
\label{MH0|0}
&   M|_{\delta=0} = O_3(\kappa), \quad
    N^0|_{\delta=0} = O_2(\kappa). &
\end{eqnarray}
The functions (\ref{MH0|0}) are analytic for $\kappa\in\calD_{\rho_\kappa}$,
$$
    \calD_{\rho_\kappa}
  = \{|\kappa_1|\le\rho_\kappa,\ldots,|\kappa_n|\le\rho_\kappa\}, \qquad
    \rho_\kappa > 0.
$$

Following the Majorant principle (Section \ref{sec:MP}), we associate with the IVP (\ref{M'H0'})--(\ref{MH0|0}) a majorant IVP:
\begin{eqnarray}
\label{M'H0'_maj}
&   \partial_\delta\bM
  = 2\bM \langle[q],\partial\rangle \bH^0, \quad
    \partial_\delta\bH^0
  = 2 \langle[q],\partial\rangle \bM e^{-2\omega_q \delta}, & \\
\label{MH0|0_maj}
&   \bM|_{\delta=0} \gg M|_{\delta=0}, \quad
    \bH^0|_{\delta=0} \gg \calH^0|_{\delta=0}. &
\end{eqnarray}

According to Lemma \ref{lem:maj1} we may choose $\bM|_{\delta=0}$ and $\bH^0|_{\delta=0}$ to be functions of $x = \kappa_1+\ldots+\kappa_n$ such that
\begin{equation}
\label{OO}
  \bM|_{\delta=0} = O(x^3)\quad\mbox{and}\quad
  \bH^0|_{\delta=0} = O(x^2).
\end{equation}
Then by Theorem \ref{theo:maj_nil} if the majorant system has a solution for all $\delta\ge 0$ then the IVP (\ref{M'H0'})--(\ref{MH0|0}) also has a solution for all $\delta\ge 0$ and
$$
  M(\kappa,\delta) \ll \bM(x,\delta), \quad
  \calH^0(\kappa,\delta) \ll \bH^0(x,\delta).
$$

Since the initial conditions (\ref{MH0|0_maj}) depend only on $x$, we may look for a solution of (\ref{M'H0'_maj})--(\ref{MH0|0_maj}) in the form $\bM=\bM(x,\delta)$, $\bH^0=\bH^0(x,\delta)$. We use the fact that for any function $F = F(x)\in\calF$ we have: $\partial_{\kappa_j} F = \partial_x F$. Therefore $\langle[q],\partial\rangle F = |q|\partial_x F$.

We put
$$
      \tau = 2 |q| \delta, \quad
       \nu = \omega_q / |q|, \quad
   \bM(x,\delta) = u(x,\tau), \quad
 \bH^0(x,\delta) = v(x,\tau).
$$
Then the IVP (\ref{M'H0'_maj})--(\ref{MH0|0_maj}) takes the form
\begin{equation}
\label{ivp_eq}
  \partial_\tau u = u\partial_x v, \quad
  \partial_\tau v = e^{-\nu\tau} \partial_x u, \qquad
  u|_{\tau=0} = O(x^3), \quad
  v|_{\tau=0} = O(x^2).
\end{equation}
Applying Theorem \ref{theo:pde} we complete the proof. \qed

\section{Technical part}

\subsection{Proof of Lemma \ref{lem:zz}}
\label{sec:zz}

Equation (\ref{z=zkap}) directly follows from Lemma \ref{lem:kqkp}. To prove (\ref{zN}), it is sufficient to take $N=\kappa^l$, where $l\in\mZ_+^n$ is arbitrary. For $\bk_q=(k_q,\overline k_q)$, $\overline k_q - k_q = q$, we have:
$$
    \{\bz^{\bk_q},\kappa^l\}
  = i\sum_{j=1}^n \big((\overline k_q)_jl_j - (k_q)_jl_j\big) \bz^{\bk_q} \kappa^{l-e_j}
  = i\bz^{\bk_q} \langle q,\partial\rangle \kappa^l.
$$

To prove (\ref{zz}), it is sufficient to consider the case $n=1$. In this case $q,p\in\mZ$ are scalar quantities.

(a) If $qp\ge0$ we have $\{\bz^{\bk_q},\bz^{\bk_p}\}=0$ and
$q\vartriangleleft p=p\vartriangleleft q=0$. Hence (\ref{zz}) holds.

(b) Suppose $p<0<-p<q$. Then
$$
  \bz^{\bk_q} = \overline z^q, \quad
  \bz^{\bk_p} = z^{-p}, \quad
  \bz^{\bk_{q+p}} = \overline z^{q+p}, \quad
  q\vartriangleleft p = 0, \quad
  p\vartriangleleft q = p.
$$
Therefore
$$
      \{\bz^{\bk_q},\bz^{\bk_p}\}
   =  - iqp\, z^{-p-1} \overline z^{q-1}
   =    i\big(q - q\vartriangleleft p - p + p\vartriangleleft q\big)\, \bz^{\bk_{q+p}} \partial\kappa^{-p} .
$$
This implies (\ref{zz}) in case (b).

(c) Suppose $p<0<q<-p$. Then
$$
  \bz^{\bk_q} = \overline z^q, \quad
  \bz^{\bk_p} = z^{-p}, \quad
  \bz^{\bk_{q+p}} = \overline z^{-q-p}, \quad
  q\vartriangleleft p = q, \quad
  p\vartriangleleft q = 0.
$$
Therefore
$$
      \{\bz^{\bk_q},\bz^{\bk_p}\}
   =  - iqp\, z^{-p-1} \overline z^{q-1}
   =    i\big(q - q\vartriangleleft p - p + [p\vartriangleleft q \big)\, \bz^{\bk_{q+p}}
           \partial\kappa^{q} .
$$
This implies (\ref{zz}) in case (c).

(d) If $p<0<-p=q$ then
$$
  \bz^{\bk_q} = \overline z^q, \quad
  \bz^{\bk_p} = z^{-p}, \quad
  \bz^{\bk_{q+p}} = 1, \quad
  q\vartriangleleft p = q/2, \quad
  p\vartriangleleft q = p/2.
$$
Therefore
$$
      \{\bz^{\bk_q},\bz^{\bk_p}\}
   =  - iqp\, z^{-p-1} \overline z^{q-1}
   =    i\big(q - q\vartriangleleft p - p + p\vartriangleleft q\big)\,
              \bz^{\bk_{q+p}} \partial\kappa^q.
$$
This implies (\ref{zz}) in case (b).

(e) The case $q<0<p$ is analogous to (b)--(d).  \qed

\subsection{Majorants}
\label{sec:maj}

For any $F,\bF\in\calF$ we say that $F\ll\bF$ iff for their Taylor coefficients we have the inequalities $|F_\bk| \le \bF_\bk$, $\bk\in\mZ_+^{2n}$.

\begin{lem}
\label{lem:maj}
Suppose $F\ll\bF$ and $\hat F\ll\hat\bF$. Then

{\bf 1}. $F+\hat F\ll\bF+\hat\bF$, $F\hat F\ll \bF\hat\bF$,
$\partial_{z_s} F\ll\partial_{z_s}\bF$, and
         $\partial_{\overline z_s} F\ll\partial_{\overline z_s}\bF$,\, $s=1,\ldots,n$.

{\bf 2}. If $F$ and $\bF$ depend on the parameter $\delta\in [\delta_1,\delta_2]$ then
$$
      \int_{\delta_1}^{\delta_2} F\, d\delta
  \ll \int_{\delta_1}^{\delta_2} \bF\, d\delta .
$$
\end{lem}

We skip an obvious proof. \qed
\smallskip

\begin{lem}
\label{lem:maj1}
Suppose $F\in\calA^\rho$, $F=O_s(\bz)$, and $\|F\|_\rho=a\rho^s$. Then
\begin{equation}
\label{geom}
  F\ll \frac{a\rho\zeta^s}{\rho - \zeta}, \qquad
  \zeta = z_1 + \ldots + z_n + \overline z_1 + \ldots + \overline z_n.
\end{equation}
\end{lem}

{\it Proof}. By Lemma \ref{lem:Hkk} we have:
$$
  F = \sum_{|\bk|\ge s} F_\bk \bz^\bk, \qquad
  |F_\bk| \le a\rho^{s - |\bk|}.
$$
Note also that $\sum_{|\bk|=j} \bz^\bk \ll \zeta^j$ for any $j\in\mZ_+$. Then
$$
     F
 \ll \sum_{|\bk|\ge s} a\rho^{s - |\bk|} \bz^\bk
 \ll a\rho^s \sum_{j=s}^\infty \frac{\zeta^j}{\rho^j}
  =  \frac{a\rho\zeta^s}{\rho - \zeta}.
$$
\qed

For any power series $F$ the corresponding majorant $\bF$ is obviously non-unique. Sometimes we will need majorant estimates of another form.

\begin{lem}
\label{lem:another_maj}
Suppose $F = \sum_{m=s}^\infty F_m x^m = O(x^s)$, $x\in\mC$ is analytic at $0\in\mC$:
$$
  \|F\|_R = c_F R^s.
$$
Then for any $\beta\in\mZ_+$
$$
  |F_m| \ll  \frac{c_F s^\beta}{m^\beta} \rho^{s-m}, \qquad
  \rho = e^{-\beta/s} R.
$$
\end{lem}

{\it Proof}. By using Lemma \ref{lem:Hkk} with $z_1=x$ and $F$ independent of
$z_2,\ldots z_n,\overline z_1,\ldots\overline z_n$, we obtain the estimate
$$
  |F_m| \le c_F R^{s-m}.
$$
It remains to use the inequality
$R^{s-m} \le s^\beta m^{-\beta} \rho^{s-m}$ for $m \ge s$.
\qed

\subsection{Majorant principle}
\label{sec:MP}

We use majorant method to obtain estimates for solutions of initial value problems (IVP) in $\calF$.

As an example consider the IVP
\begin{equation}
\label{F'=Phi}
  \partial_\delta F = \Phi(F,\delta), \qquad
  F|_{\delta=0} = \widehat F.
\end{equation}
Here $F\in\calF$ depends on the parameter $\delta$ and $\Phi$ is a map from $\calF\times\mR_+$ to $\calF$.

\begin{dfn}
\label{dfn:pr}
IVP (\ref{F'=Phi}) is said to be power regular if for any $\widehat F\in\calF$ equation (\ref{F'=Phi}) has a unique solution $F = F(\bz,\delta)\in\calF$ for all $\delta>0$
\end{dfn}

If (\ref{F'=Phi}) is a PDE then the main tool to construct a solution is the Cauchy-Kovalevskaya theorem. The system (\ref{aver1}) is not a PDE because the right-hand side includes the operators $H_\diamond\mapsto H^\sigma$, $\sigma\in\{0,+,-\}$. However in this case power regularity of the solution easily follows from nilpotent structure of the system (\ref{aver4}).

We associate with (\ref{F'=Phi}) the so-called majorant system
\begin{equation}
\label{bF'=Psi}
  \partial_\delta \bF = \Psi(\bF,\delta), \qquad
  \bF|_{\delta=0} = \widehat\bF.
\end{equation}

We put $\Phi_\bk = p_\bk\circ\Phi$ and $\Psi_\bk = p_\bk\circ\Psi$.

\begin{dfn}
\label{dfn:maj}
IVP (\ref{bF'=Psi}) is said to be a majorant IVP for (\ref{F'=Phi}) if the following two properties hold:

(a) $\widehat F\ll \widehat\bF$.

(b) For any $F\ll\bF$ and $\delta \ge 0$ we have: $\Phi_\bk(F,\delta) \ll \Psi_\bk(\bF,\delta)$ for any $\bk\in\mZ_+^{2n}$.
\end{dfn}

{\bf Majorant principle}. {\it Suppose the IVP (\ref{F'=Phi}) is power regular. Suppose also that there exists a solution $\bF = \bF(\cdot,\delta)\in\calA$ of (\ref{bF'=Psi}) on the interval
$\delta\in [0,\delta_0]$. Then (\ref{F'=Phi}) has a unique analytic solution $F$ on $[0,\delta_0]$ and $F(\cdot,\delta) \ll \bF(\cdot,\delta)$.
}
\medskip

{\bf Remarks}. 1. Definitions \ref{dfn:pr} and \ref{dfn:maj} as well as the Majorant principle obviously extend to systems of equations, where $F,\widehat F\in\calF^m$ and $\Phi:\calF^m\times\mR_+\to\calF$.

2. One may replace the first equation (\ref{bF'=Psi}) by the inequality
$\partial_\delta \bF \gg \Psi(\bF,\delta)$.

\medskip

The majorant principle presented here differs from the majorant argument used since Cauchy times. Traditionally the evolution variable (in our case $\delta$) is regarded complex as well and Taylor expansions in it are used. In our approach this variable is a real parameter in both exact solution and a majorant. Due to this we are able to obtain majorant estimates for solutions of (\ref{F'=Phi}) on large (even infinite) intervals of $\delta$.
\smallskip

\begin{theo}
\label{theo:maj_nil}
Suppose both systems (\ref{F'=Phi}) and (\ref{bF'=Psi}) have nilpotent structure. Then Majorant principle holds true.
\end{theo}

We expect that Majorant principle is valid in a much wider generality. But in this paper we are only interested in the case of systems having nilpotent form.
\smallskip

{\it Proof of Theorem \ref{theo:maj_nil}}.
Let $\bk^0$ be an index with minimal possible degree $|\bk^0|$. For example, in the system (\ref{aver4}) $|\bk^0|=3$. Nilpotent form of (\ref{F'=Phi}) implies that
$$
  0 = \partial_\delta F_{\bk^0} \ll \partial_\delta \bF_{\bk^0}.
$$
Hence $F_{\bk^0}(\delta)\ll\bF_{\bk^0}(\delta)$ for $\delta\ge 0$.

We proceed by induction in $|\bk|$. Suppose $F_{\bk}(\delta)\ll\bF_{\bk}(\delta)$, $\delta\ge 0$ provided $|\bk|< K$. For any $\bk$ such that $|\bk|=K$ we have by induction assumption and item (b) of Definition \ref{dfn:maj}:
$$
    \partial_\delta (\bF_\bk - F_\bk)
  = \Psi_\bk(\bF(\cdot,\delta),\delta) - \Phi_\bk(\bF(\cdot,\delta),\delta)
 \gg 0.
$$
Therefore
$$
    \bF_\bk(\delta)
  = \widehat\bF_\bk(\delta)
   + \int_0^\delta \Psi_\bk(\bF(\cdot,\lambda),\lambda)\, d\lambda
 \gg F_\bk(\delta)
  = \widehat F_\bk(\delta)
   + \int_0^\delta \Phi_\bk(F(\cdot,\lambda),\lambda)\, d\lambda.
$$
Here we used that arguments of $\Psi_\bk$ and $\Phi_\bk$ are known by the induction assumption.

This majorant inequality makes sense if the left-hand side is defined i.e., for any
$\delta\in [0,\delta_0]$.  \qed

\subsection{Two estimates}

We put
\begin{equation}
\label{SmN}
  S_\gamma(N) = \sum_{n_1,n_2\ge 1,\, n_1+n_2=N} \frac1{n_1^\gamma n_2^{\gamma}}.
\end{equation}

\begin{lem}
\label{lem:SmN}
Suppose $\gamma\ge 2$.
There exist a constants $\chi_\gamma > 0$ such that for any $N\ge 2$
$$
  S_\gamma(N) \le \chi_\gamma N^{-\gamma}.
$$
\end{lem}

{\it Proof}. We have:
$$
  S_\gamma(N) \le 2\sum_{1\le n_1\le n_2,\, n_1+n_2=N} \frac1{n_1^\gamma n_2^{\gamma}}
              \le 2 \Big(\frac2N\Big)^\gamma \sum_{n_1=1}^\infty \frac1{n_1^\gamma} .
$$
\qed

\begin{lem}
\label{lem:sue}
Let $\Lambda\subset\mZ_+^{2n}$ be a nonempty set, $0<\mu<\rho$, and
$F = \sum_{\bk\in\mZ_+^{2n}} F_\bk\bz^\bk \in\calA^\rho$. Consider the function
$f = \sum_{\bk\in\Lambda} F_\bk\bz^\bk$. Then
$$
  \|f\|_{\rho-\mu} \le (\rho / \mu)^{2n} \|F\|_\rho.
$$
\end{lem}

{\it Proof}. By Lemma \ref{lem:Hkk}
$$
      \|f\|_{\rho-\mu}
  \le \sum_{\bk\in\Lambda} |F_\bk| (\rho-\mu)^{|\bk|}
  \le \sum_{\bk\in\mZ_+^{2n}} \|F\|_\rho \frac{(\rho-\mu)^{|\bk|}}{\rho^{|\bk|}}
  \le \Big(\frac\rho\mu\Big)^{2n} \|F\|_\rho.
$$
\qed

\subsection{Solution of a PDE system}

Consider the PDE system (see (\ref{ivp_eq}))
\begin{equation}
\label{pde_eq1}
  \bu_\tau = \bu\bv_x, \quad  \bv_\tau = e^{-\nu\tau} \bu_x, \qquad
  \bu|_{\tau=0} = O(x^3), \quad
  \bv|_{\tau=0} = O(x^2),
\end{equation}
where the initial conditions $\bu|_{\tau=0}$ and $\bv|_{\tau=0}$ are analytic in a neighborhood of the origin.
By Lemma \ref{lem:another_maj} coefficients in the expansions
$$
  \bu(x,0) = \sum_{k=3}^\infty \bu_k(0) x^k, \quad
  \bv(x,0) = \sum_{k=2}^\infty \bv_k(0) x^k
$$
satisfy the inequalities
$$
      |\bu_k(0)|
  \le C_u k^{-3} \rho_0^{3-k}, \quad
      |\bv_k(0)|
  \le C_v k^{-4} \rho_0^{2-k},
$$
where $\rho_0>0$ may be taken arbitrarily small.

\begin{theo}
\label{theo:pde}
Suppose $\rho_0>0$ is sufficiently small. Then solution of the IVP (\ref{pde_eq1}) exists in the disk
$\Delta_r=\{x\in\mC: |x|<r\}$ for any $\tau\ge 0$, where $r>0$ is independent of $\tau$.

Moreover, there exist $Q_u,Q_v>0$ such that for any $\tau\ge 0$
$$
  \bu=O(x^3), \quad
  \bv=O(x^2), \qquad
  \|\bu\|_r \le Q_u r^3 e^{\nu\tau/2}, \quad
  \|\bv\|_r \le r^2 Q_v.
$$
\end{theo}

The proof is based on the following inductive lemma. Consider the PDE system
\begin{equation}
\label{pde_eq2}
  u_t = u v_x, \quad  v_t = \eps u_x, \qquad
  u(x,t) = \sum_{k=3}^\infty u_k(t) x^k, \quad
  v(x,t) = \sum_{k=2}^\infty v_k(t) x^k.
\end{equation}
with initial conditions
\begin{equation}
\label{ini2}
   0\le u_k(0)\le c_u k^{-3}\rho^{3-k}, \quad
   0\le v_k(0)\le c_v k^{-4}\rho^{2-k}
\end{equation}

\begin{lem}
\label{lem:induc}
Let $u,v$ be a solution of (\ref{pde_eq2}) with initial conditions, satisfying (\ref{ini2}) and let $T$ be a positive constant. Then for any $t\in [0,T]$
\begin{eqnarray}
\label{U}
     u_k(t) \le U_k(t), \quad
     U_k(t)
 &=& \frac{c_u}{k^3} \rho^{3-k} (1 + \ph e^{(\lambda + ak) t}), \qquad
     k = 3,4,\ldots, \\
\label{V}
     v_k(t) \le V_k(t), \quad
     V_k(t)
 &=& \frac{c_v}{k^4} \rho^{2-k} + \frac{c_u\psi\ph}{(k+1)^3} \rho^{2-k} e^{a(k+1) t}, \qquad
     k = 2,3,\ldots,
\end{eqnarray}
where $\ph,\psi$, and $\lambda$ do not depend on $t$ and satisfy
\begin{eqnarray}
\label{lambda}
       \ph\lambda
 &\ge& \chi\rho c_v (1 + \ph), \\
\label{ph}
       a
 &\ge& \chi\rho c_u \psi (1+\ph) e^{2aT}, \\
\label{psi}
       a \psi\ph
 &\ge& \eps (1 + \ph e^{\lambda T}).
\end{eqnarray}
The constant $\chi=\max\{\chi_2,\chi_3\}$, see Lemma \ref{lem:SmN}.
\end{lem}

We prove Lemma \ref{lem:induc} in Section \ref{sec:aux}.

\begin{cor}
\label{cor:induc}
Let $\Lambda$ be a positive constant. Then the functions $u_k$ and $v_k$ from Lemma \ref{lem:induc} satisfy the estimates
\begin{eqnarray}
\label{ini_new}
&  0\le u_k(T)\le \widehat c_u k^{-3}\widehat \rho^{3-k}, \quad
   0\le v_k(T)\le \widehat c_v k^{-4}\widehat \rho^{2-k}, & \\[.5mm]
\label{ccrho}
&  \widehat c_u = c_u e^{3(aT+\Lambda)} (1 + \ph e^{\lambda T}), \quad
   \widehat c_v = c_v + c_u e^{3aT+2\Lambda} \psi\ph / \Lambda, \quad
   \widehat\rho = \rho e^{- aT - \Lambda}. &
\end{eqnarray}
\end{cor}

Corollary \ref{cor:induc} will be proven in Section \ref{sec:aux}.
\smallskip

We choose $\tau_m = m$, $m=0,1,2,\ldots$ Let $(\bu(x,\tau),\bv(x,\tau))$ be a solution of (\ref{pde_eq1}) and $(u(x,t),v(x,t))$ a solution of (\ref{pde_eq2})$_{\eps = e^{-\nu\tau_m}}$ with initial conditions satisfying
$$
  |\bu_{k+1}(\tau_m)| \le u_{k+1}(0), \quad
  |\bv_k(\tau_m)| \le v_k(0), \qquad
  k = 2,3,\ldots
$$
Then on the interval $\tau\in [m,m+1]$
$$
  |\bu_{k+1}(\tau)| \le u_{k+1}(\tau-\tau_m), \quad
  |\bv_k(\tau)| \le v_k(\tau-\tau_m), \qquad
  k = 2,3,\ldots
$$

Hence we can obtain an estimate for $(\bu(x,\tau),\bv(x,\tau))$ applying Lemma \ref{lem:induc} and Corollary \ref{cor:induc} inductively. We will obtain sequences of positive constants
\begin{equation}
\label{sequences}
  \eps_m = e^{-\nu m}, \;\;
  T_m = 1, \;\;
  c_u(m), \;\;
  c_v(m), \;\;
  \rho_m, \;\;
  \Lambda_m, \;\;
  a_m, \;\;
  \lambda_m, \;\;
  \psi_m, \;\;
  \ph_m,
\end{equation}
where\footnote{Here we do not write the index $m$ in $u_k$ and $v_k$ for brevity.}
$$
  |\bu_k(\tau_m)| \le u_k(0) = c_u(m) k^{-3} \rho_m^{3-k}, \quad
  |\bv_k(\tau_m)| \le v_k(0) = c_v(m) k^{-4} \rho_m^{2-k}, \qquad
  m\in\mZ_+.
$$
The sequences (\ref{sequences}) have to satisfy (\ref{lambda})--(\ref{psi}) and (\ref{ccrho}).
Our aim is to prove that the sequences
$c_v(m)$ and $1 / \rho_m$ may be chosen bounded from above as $m\to\infty$ while
$c_u(m)\le Q_u\rho_m^3 e^{\nu\tau/2}$.

The quantity $\rho_0$ may be assumed to be small. We choose
$$
  a_m = \Lambda_m = 1 / (m^2 + 1), \quad
  \lambda_m = \frac{c_v(m)}{c_v(0)}, \quad
  \ph_m = \ph_0 = \sqrt{\rho_0}.
$$
Then (\ref{lambda}) holds because for sufficiently small $\rho_0$
$$
  \ph_m\lambda_m = \sqrt{\rho_0} \frac{c_v(m)}{c_v(0)} \ge \chi\rho_m c_v(m) (1 + \ph_m).
$$

Moreover, by the last equation (\ref{ccrho}) the sequence $\{\rho_m\}_{m\in\mZ_+}$ is decreasing and converges to the positive quantity $\rho_\infty = \rho_0 \exp(-\sum_{m=0}^\infty \frac2{m^2+1})$.

Conditions (\ref{ph})--(\ref{psi}) and (\ref{ccrho})\footnote
{The last equation (\ref{ccrho}) may be dropped.}
may be taken in the form
\begin{eqnarray}
\label{two}
       1
 &\ge& (m^2+1)\chi  \rho_0 c_u(m) \psi_m 2e^2, \\
\label{three_}
       \psi_m\sqrt{\rho_0}
 &=& (m^2+1) e^{-\nu m} (1 + \sqrt{\rho_0} e^{c_v(m)/c_v(0)}), \\
\label{four}
       c_u(m+1)
 &=&   c_u(m) e^{6/(m^2+1)} (1 + \sqrt{\rho_0} e^{c_v(m)/c_v(0)}), \\
\label{five}
       c_v(m+1)
 &=&   c_v(m) + (m^2+1) c_u(m) e^{5/(m^2+1)} \psi_m\sqrt{\rho_0}.
\end{eqnarray}

We define $\psi_m$ by (\ref{three_}). Then it remains to satisfy (\ref{two}), (\ref{four}), and (\ref{five}). These conditions take the form
\begin{eqnarray}
\label{two+}
       1
 &\ge& 2e^2\chi \sqrt{\rho_0} c_u(m) (m^2+1)^2 e^{-\nu m} (1 + \sqrt{\rho_0} e^{c_v(m)/c_v(0)}), \\
\label{four+}
       c_u(m+1)
 &=&   c_u(m) e^{6/(m^2+1)} (1 + \sqrt{\rho_0} e^{c_v(m)/c_v(0)}), \\
\label{five+}
       c_v(m+1)
 &=&   c_v(m) + c_u(m) e^{5/(m^2+1)} (m^2+1)^2 e^{-\nu m} (1 + \sqrt{\rho_0} e^{c_v(m)/c_v(0)}).
\end{eqnarray}

At the moment we forget about (\ref{two+}) and consider the sequences $c_u(m)$ and $c_v(m)$, satisfying (\ref{four+}), (\ref{five+}) and the initial conditions
\begin{equation}
\label{ini_cc}
  c_u(0) = C_u, \quad
  c_v(0) = C_v.
\end{equation}

Taking if necessary instead of $C_v$ a larger constant, we may assume that $C_u\le C_v$. We put
$$
  S_0 = \sum_{l=0}^\infty \frac6{l^2+1}, \quad
  S_\nu(m) = \sum_{l=0}^{m-1} (l^2 + 1)^2 e^{\nu/2 - \nu l / 2}.
$$

\begin{lem}
\label{lem:cucv}
Suppose $\rho_0>0$ is sufficiently small and $C_u\le C_v$. Then the sequences $c_u(m)$ and $c_v(m)$, computed from (\ref{four+})--(\ref{ini_cc}), satisfy the estimates
\begin{equation}
\label{cuphcv}
  c_v(m) \le C_v (1 + e^{10+S_0} S_\nu(m)), \quad
  c_u(m) \le C_u e^{S_0 + \nu m/2}.
\end{equation}
\end{lem}

{\it Proof of Lemma \ref{lem:cucv}}. Let $m_0\in\mN$ be the minimal number such that
$6 / (m_0^2 + 1) < \nu / 4$. If $\rho_0=0$ then
\begin{eqnarray*}
     c_u(m)
 &=& C_u \exp\Big(\sum_{l=1}^{m-1} \frac6{l^2+1}\Big)
 \le C_u e^{S_0}, \\
     c_v(m)
 &=& C_v + \sum_{l=1}^{m-1} c_u(l) e^{5/(l^2+1)} (l^2+1)^2 e^{-\nu l}
 \le C_v (1 + e^{5 + S_0} S_\nu(m))
\end{eqnarray*}

Therefore there exists $\rho_0>0$ such that
\smallskip

(1) inequalities (\ref{cuphcv}) hold for any $m=0,1,\ldots,m_0$,

(2) $\sqrt{\rho_0} \exp(1 + e^{10+S_0} S_\nu(\infty)) < e^{\nu/4} - 1$.
\medskip

We proceed by induction. Take any $k\ge m_0$. Suppose for any $m\le k$ inequalities (\ref{cuphcv}) hold. Then by the definition of $m_0$ and by condition (2)
$$
  e^{6/(k^2+1)} \le e^{\nu/4} \quad\mbox{and}\quad
  1 + \sqrt{\rho_0} e^{c_v(k)/c_v(0)} \le e^{\nu/4}.
$$
By using (\ref{four+}) we obtain:
$$
    c_u(k+1)
  = c_u(k) e^{6/(k^2 + 1)} (1 + \sqrt{\rho_0} e^{c_v(k)/c_v(0)})
 \le c_u(k) e^{\nu/2}.
$$
Hence the second inequality (\ref{cuphcv}) holds for $m=k+1$.

By (\ref{five+}) we have:
\begin{eqnarray*}
      c_v(k+1)
&\le& c_v(k) + c_u(k) e^{\nu/4-\nu k} (k^2+1)^2 (1 + \sqrt{\rho_0} e^{c_v(k)/c_v(0)}) \\
&\le& c_v(k) + c_u(k) (k^2+1)^2  e^{\nu/2-\nu k} \\
&\le& C_v(1 + e^{10 + S_0} S_\nu(k)) + C_u e^{S_0+\nu k/2} (k^2+1)^2  e^{\nu/2-\nu k} \\
&\le& C_v(1 + e^{10 + S_0} S_\nu(k+1)).
\end{eqnarray*}
Hence the first inequality (\ref{cuphcv}) holds for $m=k+1$. \qed

Finally turn to inequality (\ref{two+}). If $\rho_0$ is taken satisfying condition (2) in the proof of Lemma \ref{lem:cucv} then by Lemma \ref{lem:cucv} we have:
$1 + \sqrt{\rho_0} e^{c_v(m) / c_v(0)} \le e^{\nu/4}$ and $c_u(m)$ satisfies (\ref{cuphcv}). Then (\ref{two+}) takes the form
$$
  1 \ge 2e^2\chi \sqrt{\rho_0} C_u e^{S_0} (m^2+1)^2 e^{\nu/2 - \nu m/2}.
$$
Taking if necessary a smaller $\rho_0$, we can satisfy this inequality for any $m\ge 0$. This finishes proof of Theorem \ref{theo:pde}.

\subsection{Auxiliary IVP}
\label{sec:aux}

In this section we prove Lemma \ref{lem:induc} and Corollary \ref{cor:induc}.
\smallskip

{\bf Proof of Lemma \ref{lem:induc}}. Logic of our argument is as follows:
\smallskip

(a) we prove (\ref{U})$_{k=3}$,

(b) we derive (\ref{V})$_k$ from (\ref{U})$_{k+1}$,

(c) we derive (\ref{U})$_k$ from (\ref{V})$_{n\le k-2}$ and (\ref{U})$_{n\le k-1}$.
\medskip

The system (\ref{pde_eq2}) is equivalent to the following infinite ODE system:
\begin{equation}
\label{pde->ode}
  (u_k)_t = \sum_{\alpha+\beta=k+1} \beta u_\alpha v_\beta, \quad
  (v_k)_t = \eps (k+1) u_{k+1}.
\end{equation}

(a) The equation $(u_3)_t=0$ and the first estimate (\ref{ini2})$_{k=3}$ imply (\ref{U})$_{k=3}$.

(b) By (\ref{ini2}) estimate (\ref{V})$_k$ holds for $t=0$. To prove (\ref{V})$_k$ for arbitrary $t\in [0,T]$, we assume that (\ref{U})$_{k+1}$ is valid. Then it is sufficient to prove that $(v_k)_t\le (V_k)_t$. By (\ref{pde->ode}) we have:
\begin{eqnarray*}
        (v_k)_t
  &\le& \frac{\eps c_u}{(k+1)^2} \rho^{2-k} (1 + \ph e^{(\lambda+a(k+1)) t}), \\
        (V_k)_t
   &=&  \frac{a c_u\psi\ph}{(k+1)^2} \rho^{2-k} e^{a(k+1) t}.
\end{eqnarray*}
Then the inequality $(v_k)_t\le (V_k)_t$ follows from (\ref{psi}).

(c) By (\ref{ini2}) estimate (\ref{U})$_k$ holds for $t=0$. To prove (\ref{U})$_k$ for arbitrary $t\in [0,T]$, we assume that (\ref{V})$_{n\le k-2}$ and (\ref{U})$_{n\le k-1}$ are valid. It is sufficient to prove that $(u_k)_t\le (U_k)_t$. By (\ref{pde->ode})
\begin{eqnarray*}
        (u_k)_t
  &\le& \sum_{\alpha+\beta=k+1} \rho^{4-k} \frac{c_u}{\alpha^3}
              \Big( 1 + \ph e^{(\lambda+a\alpha) t} \Big)
              \Big( \frac{c_v}{\beta^3} + \frac{\beta c_u \psi\ph}{(\beta+1)^3} e^{a(\beta+1) t} \Big) \\
   &=&  \rho^{4-k} (A_1 + A_2)  , \\
        A_1
   &=&  c_u c_v \sum_{\alpha+\beta=k+1} \big( 1 + \ph e^{(\lambda+a\alpha)t} \big)
                                               \frac{1}{\alpha^3 \beta^3},  \\
        A_2
   &=&  c_u^2 \ph\psi \sum_{\alpha+\beta=k+1}
           \Big( 1 + \ph e^{(\lambda+a\alpha) t} \Big) \frac{e^{a(\beta+1) t}}{\alpha^3 (\beta+1)^2}.
\end{eqnarray*}
We estimate $A_1$ and $A_2$ separately. By using Lemma \ref{lem:SmN} we have:
\begin{eqnarray*}
        A_1
 &\le&  c_u c_v \big( 1 + \ph e^{(\lambda+ak) t} \big)
           \sum_{\alpha+\beta=k+1} \frac{1}{\alpha^3 \beta^3}
 \;\le\; \frac{c_u c_v (1 + \ph) e^{(\lambda+ak) t}}{(k+1)^3} \chi ,  \\
        A_2
 &\le&  c_u^2 \ph\psi \sum_{\alpha+\beta=k+1}
           ( 1 + \ph ) \frac{e^{(\lambda + a(k+2)) t}}{\alpha^3 (\beta+1)^2}
 \;\le\;  \frac{c_u^2 \ph\psi(1 + \ph)}{(k+1)^2} e^{(\lambda + a(k+2)) t}  \chi
\end{eqnarray*}

We also obtain from (\ref{U})
$$
      (U_k)_t
   =  \frac{c_u (\lambda+ak)}{k^3} \rho^{3-k} \ph e^{(\lambda+ak) t}.
$$

Hence the inequality $(u_k)_t\le (U_k)_t$ follows from two inequalities
$$
  \rho^{4-k} A_1 \le \frac{c_u \lambda}{k^3} \rho^{3-k} \ph e^{(\lambda+ak) t} \quad
  \mbox{and}\quad
  \rho^{4-k} A_2 \le \frac{c_u a}{k^2} \rho^{3-k} \ph e^{(\lambda+ak) t}.
$$
The first one follows from (\ref{lambda}) while the second from (\ref{ph}). \qed
\smallskip

{\bf Proof of Corollary \ref{cor:induc}}.

To prove the first estimate (\ref{ini_new}) we have to check that
$$
      \frac{c_u}{k^3} \rho^{3-k} \big(1 + \ph e^{(\lambda + ak)T}\big)
  \le \frac{\widehat c_u}{k^3} \widehat\rho^{3-k}.
$$
By (\ref{ccrho}) this inequality is equivalent to
$$
      1 + \ph e^{(\lambda + ak)T}
  \le e^{3(aT+\Lambda)} (1 + \ph e^{\lambda T}) e^{-(3-k)(aT+\Lambda)}.
$$
which obviously holds.

To prove the second estimate (\ref{ini_new}) we check that
$$
       \frac{c_v}{k^4} \rho^{2-k}
      + \frac{c_u\psi\ph}{(k+1)^3} \rho^{2-k} e^{a(k+1)T}
  \le \frac{\widehat c_v}{k^4} \widehat\rho^{2-k}.
$$
By (\ref{ccrho}) this inequality follows from
$$
      \frac{c_u\psi\ph}{k^3} \rho^{2-k} e^{a(k+1)T}
  \le \frac{c_u e^{3aT+2\Lambda}\psi\ph}{\Lambda k^4} \rho^{2-k} e^{-(2-k)(aT+\Lambda)}
$$
which is equivalent to $k\Lambda\le e^{k\Lambda}$.  \qed

\section{Plans}

\begin{itemize}

\item Take instead of $\xi$ the operator $\xi_M=p_{|\bk|\le M}\circ\xi$. Do we normalize everything for $|\bk|\le M$?

\item  Normalization outside $S_{-M,M}$.

\item Gevrey estimates for $H^0(+\infty)$, $\calH^Q(\cdot,\tau)$.

\item  $\phi^\tau(\calA^\rho) \subset \calA^{f(\rho,\tau)}$. Estimate $f$ better than
$f\sim 1/\delta$. Does $H$ approach the normal form in $D_{f(\rho,\tau)}$ ?

\item  Is $W_0$ somewhat specific?

\item  Group of formal canonical transformations, its action on $\calF$, on $W_N$.

\item  $\widehat H\in\calA$, one DOF. Do we have convergence?

\item  $\calN$ is asymptotically stable. Do we have a Lyapunov function?

\item  Symmetries

\item  Example: $\widehat H\in\calA$, $BNF(\widehat H)\not\in\calA$.

\item  Normalization of invariant (may be, asymptotic) manifolds.

\item Theorem.  Suppose $\|\widehat H\|_{\rho_0}=c_H$. Then there exists $\widehat G = O_N(\bz)$,
      $\|\widehat G\|_\rho \le C(c_H,N)$ such that the system with Hamiltonian
      $\widehat H+\widehat G$ is completely integrable.
\end{itemize}

\end{document}